\documentclass[letterpaper, 10 pt, conference]{ieeeconf}  

\IEEEoverridecommandlockouts                              

\renewcommand{\QED}{\QEDopen}

\usepackage{graphics} 
\usepackage{epsfig} 
\usepackage{amsmath} 
\usepackage{amssymb}  
\usepackage{color}

\usepackage{epstopdf}
 
\title{\LARGE \bf
	Adaptive sliding mode control without knowledge of uncertainty bounds}

\author{Yi-Wen Liao $^{1}$, Selina Pan $^{2}$, Francesco Borrelli$^{1}$ and J. Karl Hedrick $^{1}$
	\thanks{$^{1}$Yi-Wen Liao, Francesco Borrelli and J. Karl Hedrick are with the Department of Mechanical Engineering, University of California, Berkeley, CA 94720, USA. Email: {\tt\small \{ywliao, fborrelli, khedrick\}@berkeley.edu}.}
	\thanks{$^{2}$ Selina Pan is with Research and Innovation CENTER, Ford Motor Company, Palo Alto, CA 94304, USA. Email: {\tt\small span6@ford.com}}
}

\begin{document}

	\maketitle
	\thispagestyle{empty}
	\pagestyle{empty}

	\begin{abstract}
		This paper proposes a new adaptation methodology to find the control inputs for a class of nonlinear systems with time-varying bounded uncertainties. The proposed method does not require any prior knowledge of the uncertainties including their bounds. 
		The main idea is developed under the structure of adaptive sliding mode control; an update law decreases the gain inside and increases the gain outside a vicinity of the sliding surface.
		The semi-global stability of the closed-loop system and the adaptation error are guaranteed by Lyapunov theory. 
		The simulation results show that the proposed adaptation methodology can reduce the magnitude of the controller gain to the minimum possible value and smooth out the chattering.      
	\end{abstract}   
	\section{INTRODUCTION}
	Sliding mode control \cite{edwards1998sliding}-\cite{utkin2009sliding} has been recognized as one of the effective nonlinear control methods due to its robustness to uncertainties and its guarantee of finite time convergence. However, the design procedure requires the knowledge of the bound on the uncertainties, which, from a practical point of view, is usually hard to acquire. This results in an uncertainty bound that is often overestimated and hence leads to an undesirable large control gain in the discontinuous sliding term. Consequently, the system will suffer from large magnitude chattering behaviors \cite{boiko2005analysis}.
	
	To reduce this kind of ``zig-zag" motion, several methods have been proposed, which include the boundary layer technique \cite{yao1996smooth} and the ``equivalent" control method \cite{utkin1992sliding}-\cite{tseng2010chattering}. The first, proposed by Yao and Tomizuka, approximates the discontinuous signum function by a high-slope saturation function with a desired thickness of the boundary layer. The second, shown by Utkin as well as Tseng and Chen, replaces the discontinuous signum function with a low-pass filter. Although we can get a continuous sliding controller from these methods, the guarantee of global asymptotic stability is sacrificed \cite{slotine1991applied}. In addition, both of these approaches require prior knowledge of the bound on the uncertainties. To avoid this, we can make use of the adaptive control strategy \cite{sastry2011adaptive}-\cite{krstic1995nonlinear} to estimate the unknown parameters.  
	Common methods of estimation include recursive least squares and gradient descent. A more direct way is to derive the update laws from Lyapunov stability theory and analyze the convergence performance. The update laws will use the current information to modify the control input in real time.
	Because of the advantage of not overestimating the bound on the uncertainties, many adaptation approaches combined with sliding mode control have been developed to tune the sliding gains.
	The adaptation law proposed in \cite{huang2008adaptive} is proportional to the tracking error. It shows that the system will converge to the sliding surface within a finite time. However, the sliding gain will approach infinity since the ideal sliding mode does not exist. 
	In \cite{hussain2004adaptive}, neural networks model the uncertainties of the system and the resulting controller is implemented on a two-tank level control system. The results show that it can enable a lower switching gain and eliminate the chattering with a thin boundary layer. However, it requires an off-line training process and cannot guarantee stability.
	Another gain-adaptation algorithm is proposed by using a sliding mode disturbance observer \cite{hall2006sliding}, but it has the drawback of requiring the knowledge of uncertainty bounds.
	The objective of this paper is to provide an adaptive control methodology for a class of nonlinear systems with uncertainties. Note that the uncertainties should be bounded but the prior knowledge of the bound is unknown.
	
	The paper is organized as follows. 
	In Section II, we review two adaptive sliding mode control strategies proposed in \cite{utkin2013adaptive}.  
	Then, the new adaptive sliding mode control is introduced in Section III. The stability analysis is provided.
	Finally, the performance of the proposed method is demonstrated with two examples and compared with one of the existing methods.
	
	\section{PRELIMINARIES}
	\subsection{Problem Statement}
	Consider a nonlinear system given by:
	\begin{align}
	\begin{cases}
	\dot{\mathbf{x}}(t)= \mathbf{f}(\mathbf{x},t)+\mathbf{l}(\mathbf{x},t)u(t) \\
	y(t)= \mathbf{c}(\mathbf{x},t)
	\end{cases}~
	\mathbf{x}(0) = \mathbf{x_0}, ~t \ge 0 \label{Problem}
	\end{align}
	where $\mathbf{x}(t) = [x_1(t), x_2(t), ..., x_n(t)]^\text{T}\in \mathcal{X} \subset \mathbb{R}^n$ is the state vector, $u(t) \in \mathbb{R}$ is the control input and $y(t) \in \mathbb{R}$ is the system output. $\mathbf{f}(\mathbf{x},t)$ and $\mathbf{l}(\mathbf{x},t)$ are bounded and sufficiently smooth functions which describe the model of the system. Assume that both of them contain unmeasured model uncertainties which satisfy the ``matching condition" for conventional sliding mode control \cite{dravzenovic1969invariance}. Additionally, to guarantee the controllability, $\mathbf{l}(\mathbf{x},t)$ should be $\neq \mathbf{0}$ for all $(\mathbf{x},t)\in \mathcal{X}\times\mathbb{R}^+$. 
	
	The common goal of the control problem is to guide the output $y(t)$ along a desired trajectory, $y_d(t)$, or around the origin. To design the sliding mode control, first we define a stable sliding surface $s(\mathbf{x},t)$~\cite{utkin2009sliding} with a relative degree equal to $1$ with respect to $u(t)$. Then, we obtain the time derivative of $s(\mathbf{x},t)$ as 
	\begin{align}
	\dot{s}(\mathbf{x},t) & = \frac{\partial s(\mathbf{x},t)^T}{\partial \mathbf{x}}\dot{x}+\frac{\partial s(\mathbf{x},t)}{\partial t} \nonumber\\
	&= h(\mathbf{x},t)+g(\mathbf{x},t)u(t) \label{eq:sy}
	\end{align}
	where 
	\begin{align}
	    h(\mathbf{x},t) &=  \frac{\partial s(\mathbf{x},t)}{\partial t}+\frac{\partial s(\mathbf{x},t)^T}{\partial \mathbf{x}}\mathbf{f}(\mathbf{x},t)\nonumber\\
	    g(\mathbf{x},t) &=  \frac{\partial s(\mathbf{x},t)^T}{\partial \mathbf{x}}\mathbf{l}(\mathbf{x},t). \nonumber
	\end{align}
	\begin{figure}[t]
		\centering
		\includegraphics[width = 2.6in]{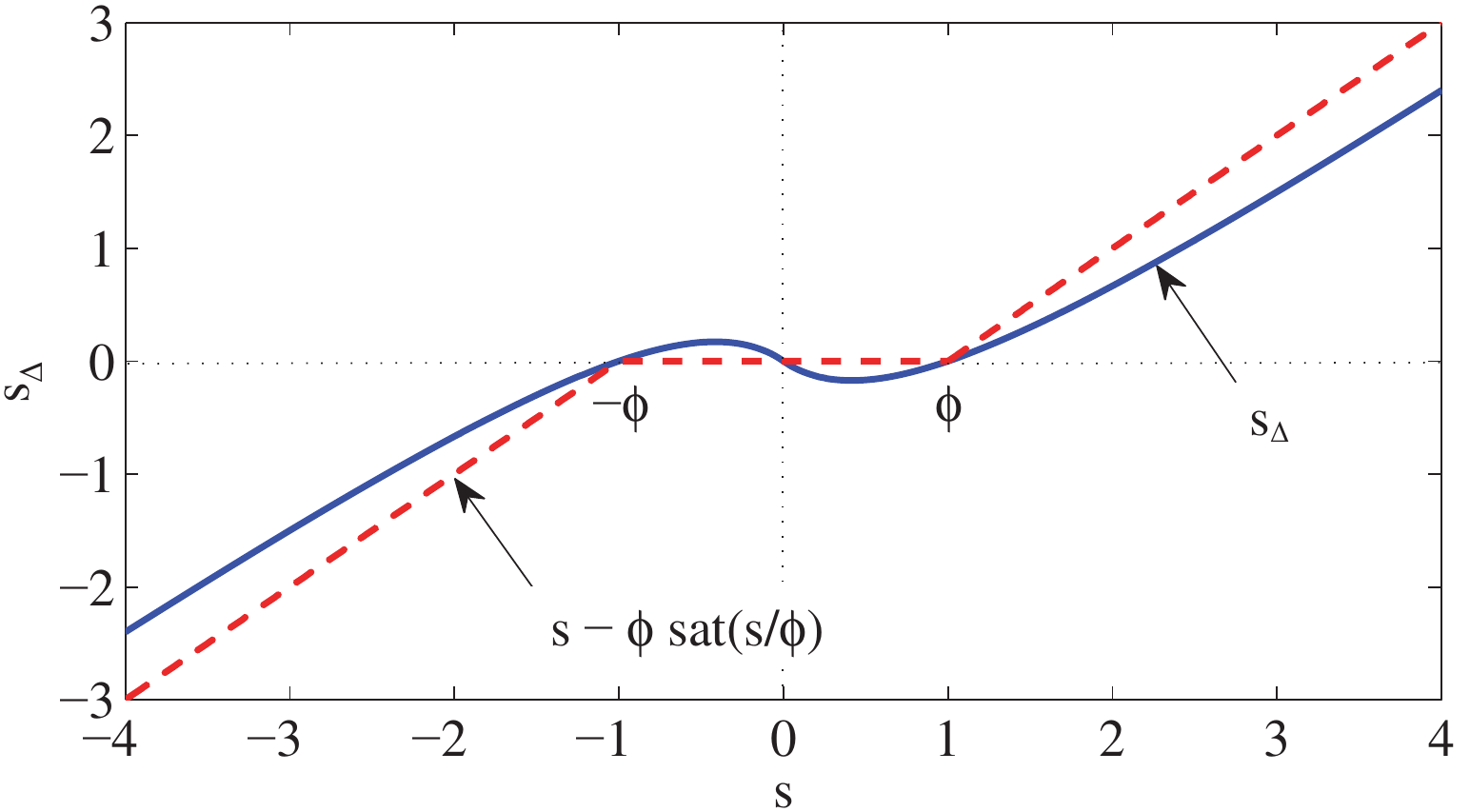}
		\caption{Comparison between the standard delta function with linear saturation function $\text{sat}(s/\phi)$ and the proposed delta function $s_{\Delta}$.}			\label{fig:ds}
	\end{figure}	
	To handle the modeling uncertainties and unknown disturbances, we rewrite the model \eqref{eq:sy} with an addictive time-varying function, $\Delta f(\mathbf{x},t)$:
	\begin{align}
	\dot{s}(\mathbf{x},t) &= h(\mathbf{x},t)+g(\mathbf{x},t)u(t)+\Delta f(\mathbf{x},t).\nonumber
	\end{align}
	The term of $\Delta f(\mathbf{x},t)$ represents the overall uncertainty of the system and satisfy the following inequality:
	\begin{align}
	|\Delta f(\mathbf{x},t)| \leq \varsigma(\mathbf{x},t) \leq \mu\nonumber
	\end{align}
	where $\mu$ is the \textit{unknown} upperbound.

	The objective in this paper is to design a control law which can adapt the time-varying uncertainty, $\varsigma(\mathbf{x},t)$, in order to reduce the chattering behavior in conventional sliding mode control, but still preserve its own strength in the guarantee of robustness and stability. 
	\subsection{Adaptive Sliding Mode Control Revisit}
	As is common for sliding mode control, the controller is designed as 
	\begin{align}
	u = -K\text{sgn}(s) \label{eq:SMC control}
	\end{align}
	where the controller gain, $K$, is the design parameter which should be greater than or equal to the uncertainty bound, $\mu$.
	\begin{align}
	\text{sgn}(s) \doteq \begin{cases}
	1~~~~~~~~~~~~~~\text{if}~~~s>0\\
	-1~~~~~~~~~~~~\text{if}~~~s<0\\
	0~~~~~~~~~~~~~~\text{if}~~~s=0
	\end{cases} \nonumber
	\end{align}
	is the discontinuous switching function \cite{khalil2002nonlinear}. 
	As mentioned in the introduction, having a poor estimation on the upperbound, $\mu$, will lead to a larger chattering behavior in the system response.
    \begin{figure}[t]
    	\centering
    	\includegraphics[width = 2.6in]{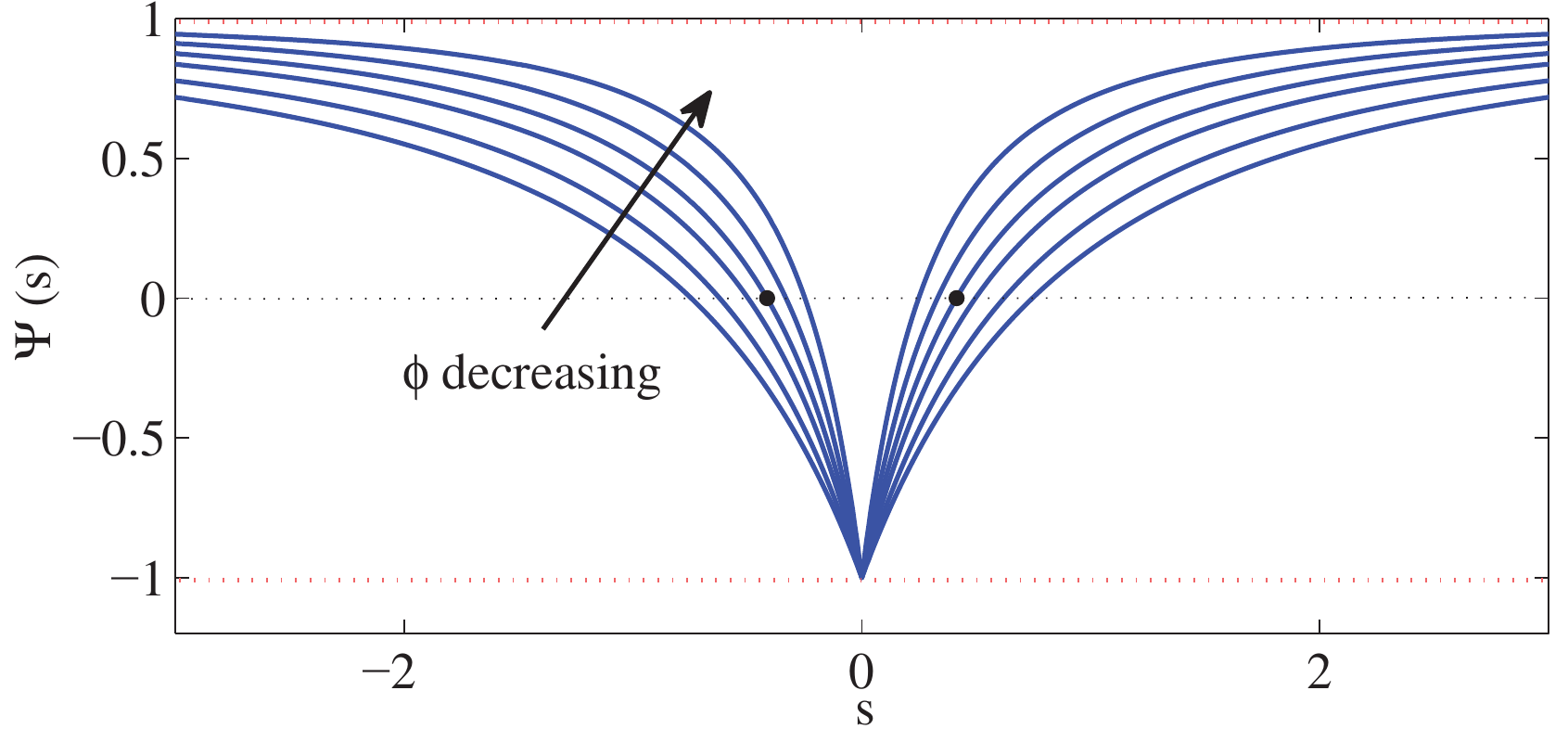}
    	\caption{The plot of $\Psi(s)$ with varying thicknesses of the boundary layer $\phi$. The dotted red lines indicate the upper and lower bounds of $1$ and $-1$.}
    	\label{fig:delta_s}
    \end{figure}
	Thus, the main goal of the adaptive sliding mode control is to reduce the magnitude of the controller gain, $K$, to its minimum admissible value. In other words, the controller gain is not a constant anymore; instead, it can be tuned and modified with time.
	The method proposed in \cite{utkin2013adaptive} is based on the use of ``equivalent" control: once sliding mode occurs, the uncertainty magnitude can be evaluated and adequately tuned by a low-pass filter:
	\begin{align}
		[\text{sgn}(s)]_{\text{eq}}\doteq z\in (0,1):~~\tau \dot{z}+z = \text{sgn}(s(\mathbf{x},t)). \label{eq:equivalentControl}
	\end{align}
	To preserve sliding mode and minimize the chattering, the ideal gain $K(t)$ should tend to $\Delta f(t)$ and be slightly greater than $\Delta f(t)$. So, the design idea of the adaptation would be:
	\begin{align}
	K(t) \approx |\Delta f(t)|/\alpha,~~\alpha \in (0,1)\nonumber
	\end{align}
	where $\alpha$ is very close to $1$. According to this, the minimal possible value of the gain $K$ can be found using the following adaptation algorithm:
	\begin{align}
	\dot{K} &= \nu K\text{sgn}(\delta)-M[K-K^+]_++M[\epsilon-K]_+ \label{eq:adaptationLaw}
	\end{align}
	with
	\begin{align}
	\delta &\doteq \left|[\text{sgn}(s(\mathbf{x},t))]_\text{eq}\right|-\alpha,~~~\alpha \in (0,1) \nonumber\\
	[z]_+ &\doteq \begin{cases}
	1 ~~\text{if}~~z\geq 0 \\
	0 ~~\text{if}~~z<0,
	\end{cases}~M>\nu K^+,~K^+ > \mu ,\nu >0. \nonumber
	\end{align}
	$\epsilon>0$ is a preselected minimal value of $K$ and $K^+$ is the uncertainty bound. Once sliding mode with respect to $s(\mathbf{x},t)$ is established, the adaptation law (\ref{eq:adaptationLaw}) allows the gain $K$ to vary in the range of $[\epsilon, K^+]$ and to be slightly greater than the current uncertainty $\Delta f(t)$. This guarantees an ideal sliding motion. 
	
	Another strategy is proposed in \cite{plestan2010new} which is similar to what we have just introduced above. Instead of using the ``equivalent" control method to estimate the boundary of the uncertainties, consider the adaptation law:
	\begin{align}
	\dot{K} = \begin{cases}
	\bar{K}|s(\mathbf{x},t)|\text{sgn}(|s(\mathbf{x},t)|-\epsilon)~&\text{if}~K>\kappa\\
	0 ~~~~~~~~~~~~~~~~~~~~~~~~~~~~~&\text{if}~K \le\kappa 
	\end{cases}\label{eq:Adap}
	\end{align}
	with $\bar{K}>0$, $\epsilon>0$ and a small enough value of $\kappa>0$ that ensures a positive value of $K$. 
	According to (\ref{eq:Adap}), $K$ will decrease if $|s(\mathbf{x},t)|<\epsilon$. In other words, the gain $K$ will be kept at the smallest level that allows a given certain amount of accuracy which means we can only guarantee semi-global stability of the system. However, the big advantage of this method is that it does not require the knowledge of the uncertainty bound.
	\begin{figure}[t]
		\centering
		\includegraphics[width = 2.4in]{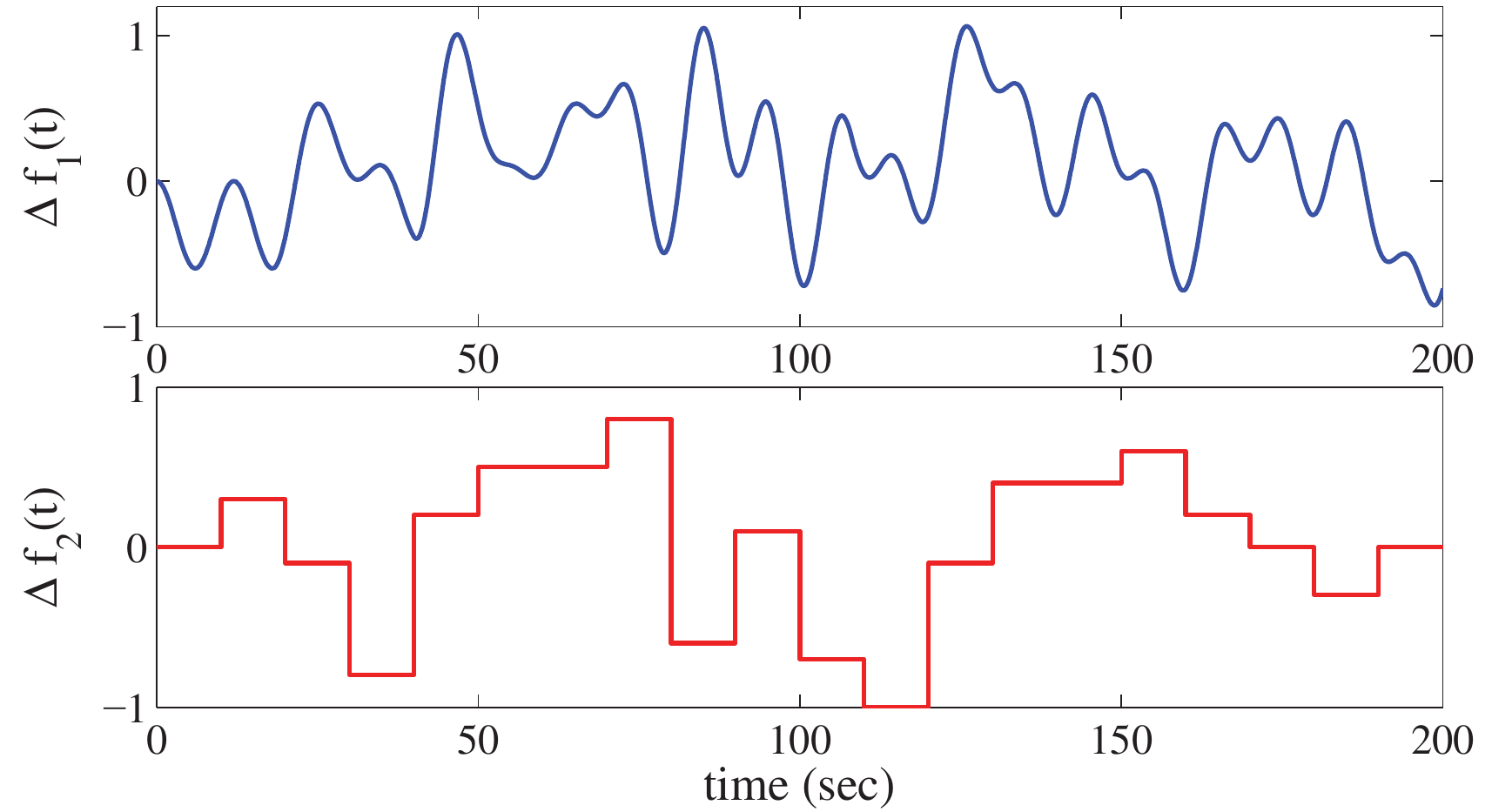}
		\caption{The uncertainty functions $\Delta f_1(t)$ and $\Delta f_2(t)$ vs. time}
		\label{fig:df}
	\end{figure}  
	\begin{figure}[t]
		\centering
		\includegraphics[width = 2.8in]{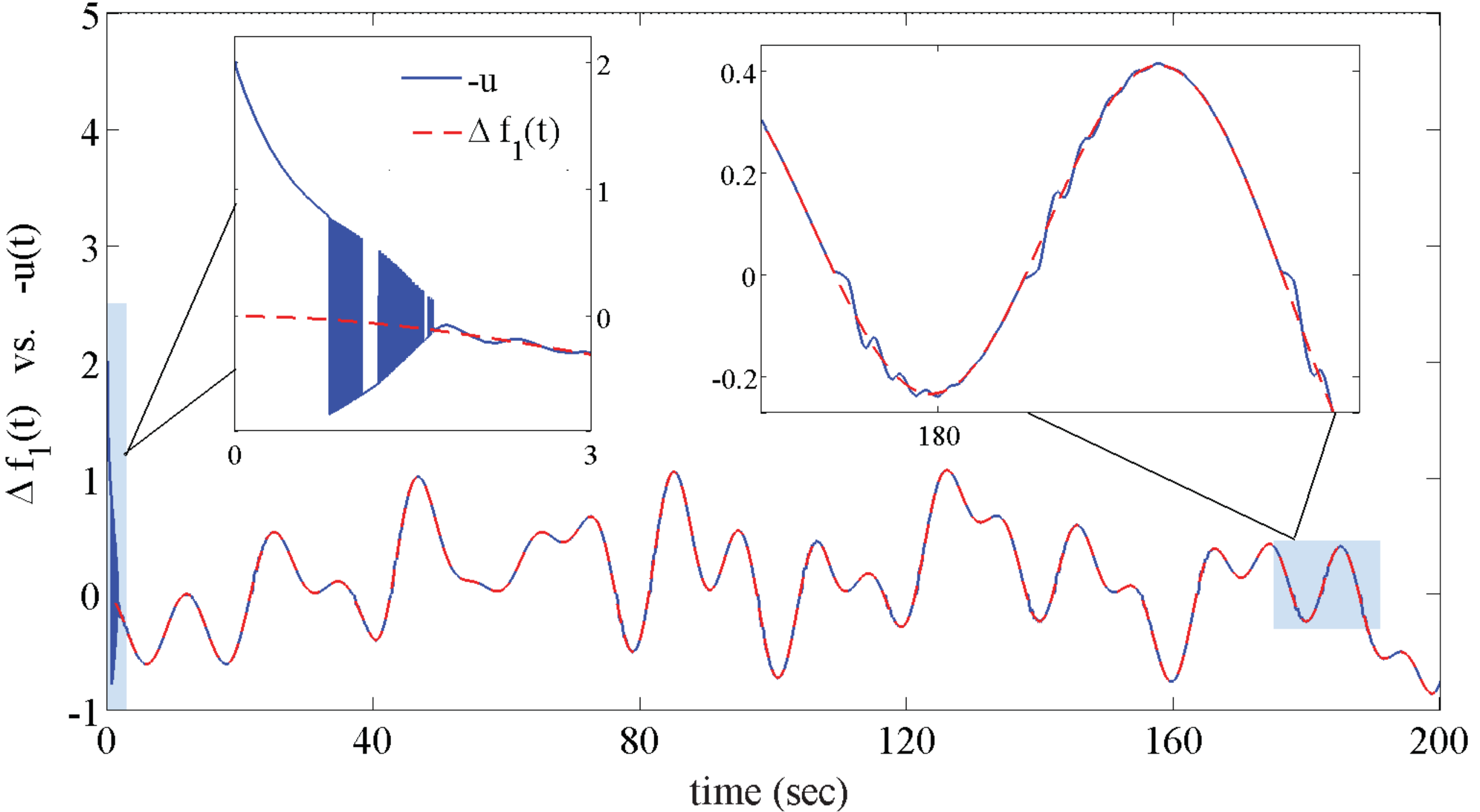}
		\caption{The adaptation performance for the smooth uncertainties $\Delta f_1(t)$}
		\label{fig:f1}
	\end{figure}
	\begin{figure}[t]
		\centering
		\includegraphics[width = 2.4in]{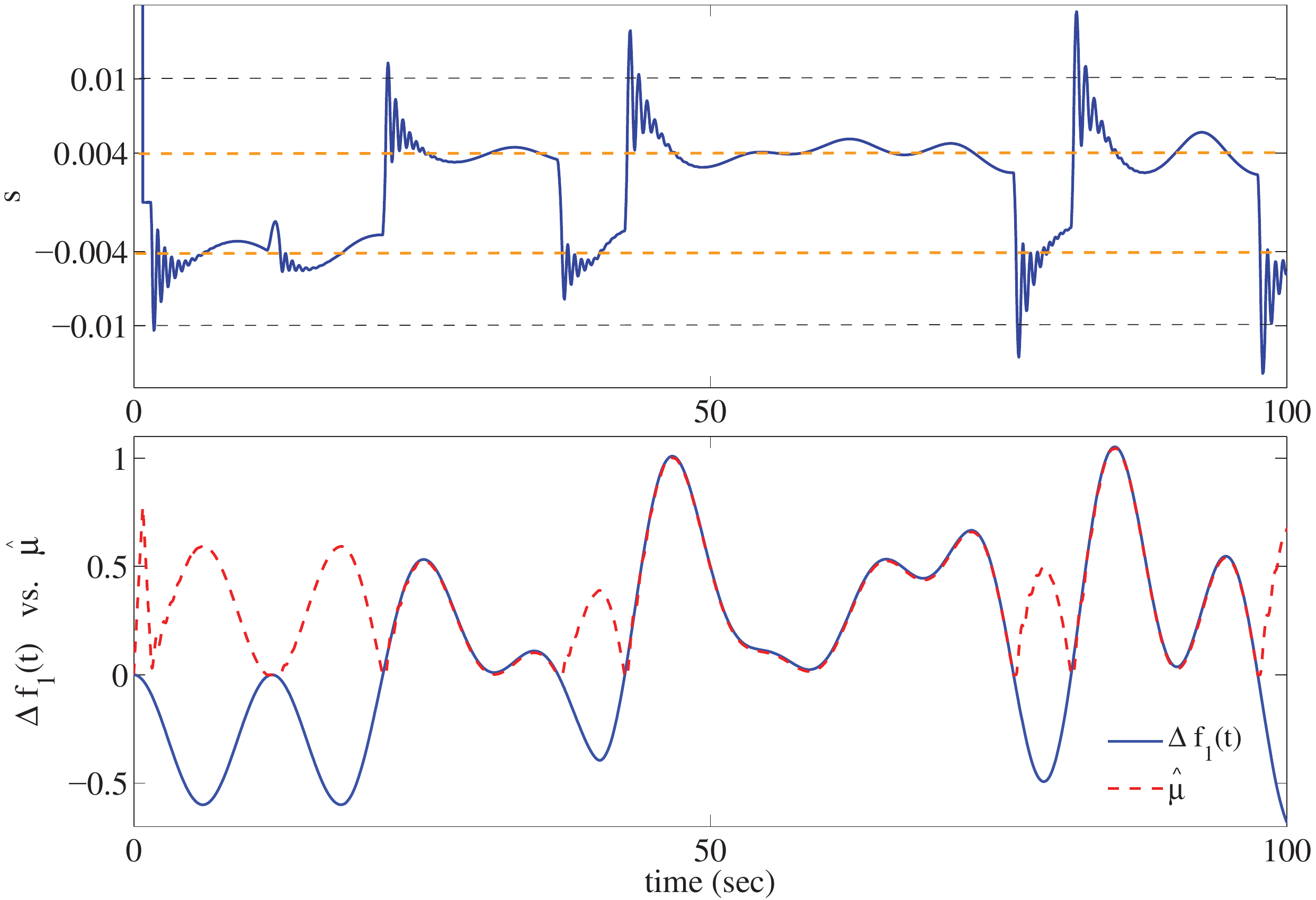}
		\caption{The convergence value of the sliding variables $s(t)$ vs. the adaptation gain $\hat{\mu}$ for the smooth uncertainties $\Delta f_1(t)$}
		\label{fig:f2}
	\end{figure}
	\begin{figure}[t]
		\centering
		\includegraphics[width = 2.5in]{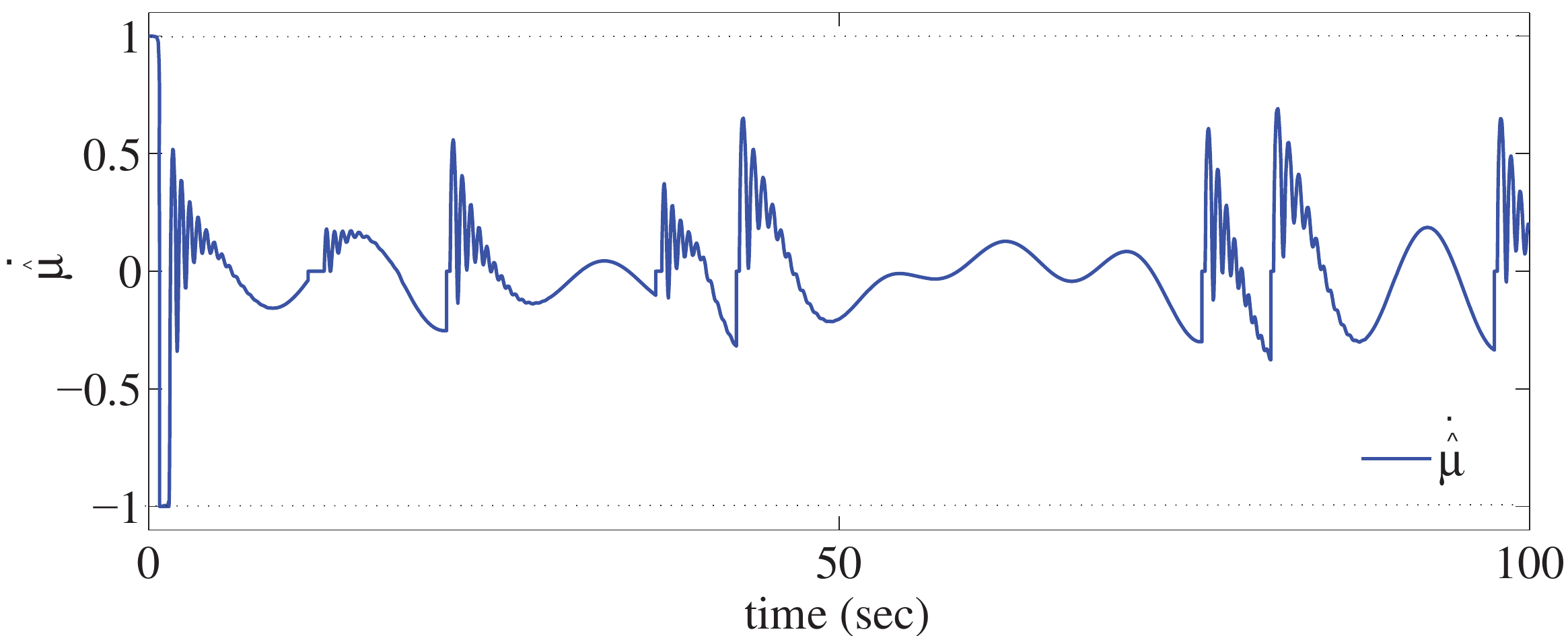}
		\caption{The changing rate of the adaptation gain, $\dot{\hat{\mu}}$ for the smooth uncertainties $\Delta f_1(t)$.}
		\label{fig:f3}
	\end{figure}
	\section{NEW ADAPTIVE SLIDING MODE CONTROL}
	\subsection{Motivation and Design Idea}
	Although the method proposed in \cite{plestan2010new} has a big advantage for not requiring the knowledge of the uncertainty bound, the adaptation algorithm (\ref{eq:Adap}) will introduce a discontinuous changing rate in control $u(t)$ at $|s(\mathbf{x},t)|=\epsilon$, which is not realistic for many actuators. Another problem is that the adaptation law is a linear function of $|s(\mathbf{x},t)|$, which makes the adaptation rate (increasing and decreasing) not quite even. In other words, the adaptive rate will become unreasonably large or too small when $s(\mathbf{x},t)$ is in the reaching phase or converging around zero. As a result, the common problem of the sliding mode control, chattering behaviors, will be easily exhibited in the steady state. To overcome this, we propose another methodology incorporated with the concept of the boundary layer.
	Using the similar idea that the gain will increase outside and decrease inside the small region around the sliding surface, we introduce a special delta function which ``roughly" denotes the distance of $s$ from the boundary layer. The function is defined as:
	\begin{align}
	s_{\Delta}(\mathbf{x},t) \doteq s - \frac{2s\phi}{|s|+\phi}
	\end{align} 
    where $\phi >0$ is a design parameter indicating the thickness of the boundary layer. 
    
    It is worth noting that, in comparison with the classical delta function defined by the saturation function $\text{sat}(s/\phi)$, the new delta function has a similar shape but with nonzero values inside the boundary layer (See Fig. \ref{fig:ds}). 
    There are three main advantages of using $s_{\Delta}(\mathbf{x},t)$ to derive the adaptation law. First, instead of blindly tuning a time constant $\tau$ of the low-pass filter in (\ref{eq:equivalentControl}) or the adaptation gain $\bar{K}$ in (\ref{eq:Adap}), the new adaptation law provides a smooth adaptation process based on the feedback information from $s_{\Delta}(\mathbf{x},t)$. 
    Second, unlike the chattering behavior in many adaptive sliding mode control algorithms, it can alleviate the chattering with a simple parameter tuning method.
    Finally, the stability proof can be done in a clean and relatively easy way.    
	\subsection{New Adaptation Control Law}
	Consider the same problem described in Section II.A with the sliding surface $s(\mathbf{x},t)$ defined in the same way as listed in (\ref{eq:sy}). The following theorems describe the stability property with the adaptation law based on the delta function we proposed.	
	
	\textbf{Theorem III.1:} Given the system (\ref{Problem}) implemented with the following feedback control and adaptive update laws, the closed-loop state $s$ will approach the boundaries of the domain $\mathcal{S}=\left\{s\in \mathbb{R}\text{, } |s| \ge \eta\right\}$ for $\eta = (\sqrt{2}-1)\phi$.
	\begin{eqnarray}
	\begin{aligned}
		~u &= -\frac{1}{g(\mathbf{x},t)}\left[\hspace{0.3mm}h(\mathbf{x},t)+ks+\hat{\mu}\text{sgn}(s)\hspace{0.3mm}\right] \\
		~\dot{\hat{\mu}} &= \begin{cases}
		 \frac{1}{\rho}\left[1-\frac{2 \phi^2}{(|s|+\phi)^2}\right]\text{~~if~~} \hat{\mu}\ge 0\\
		 0\text{~~~~~~~~~~~~~~~~~~~~if~~} \hat{\mu} < 0 
		\end{cases} \hat{\mu}(0) = \hat{\mu}_0
	\end{aligned}\label{eq:control}
	\end{eqnarray}
	where $\rho > 0$ is the adaptation gain, $k > 0$ is the feedback control gain and $\hat{\mu}_0 > 0$ is the initial guess of the sliding gain.\\
	\textit{Proof:}              
	We first calculate the time derivative of the sliding surface $s$ and $s_{\Delta}$ from (\ref{fig:ds}) and (\ref{eq:control}) as follows:
	\begin{align}
	\dot{s} &= \Delta f - \hat{\mu}\text{sgn}(s)-ks \nonumber\\
	\dot{s}_{\Delta} & = \dot{s}\left[1-\frac{2\phi^2}{(|s|+\phi)^2}\right]. \nonumber
	\end{align}	
	Consider the following Lyapunov function candidate:
	\begin{align}
	V(s,\hat{\mu}) =\text{sgn}(s) s_{\Delta}+\frac{1}{2}\rho(\mu-\hat{\mu})^2. \label{eqn:Lyapunov} 
	\end{align}
	\begin{figure}[t]
		\centering
		\includegraphics[width = 2.85in]{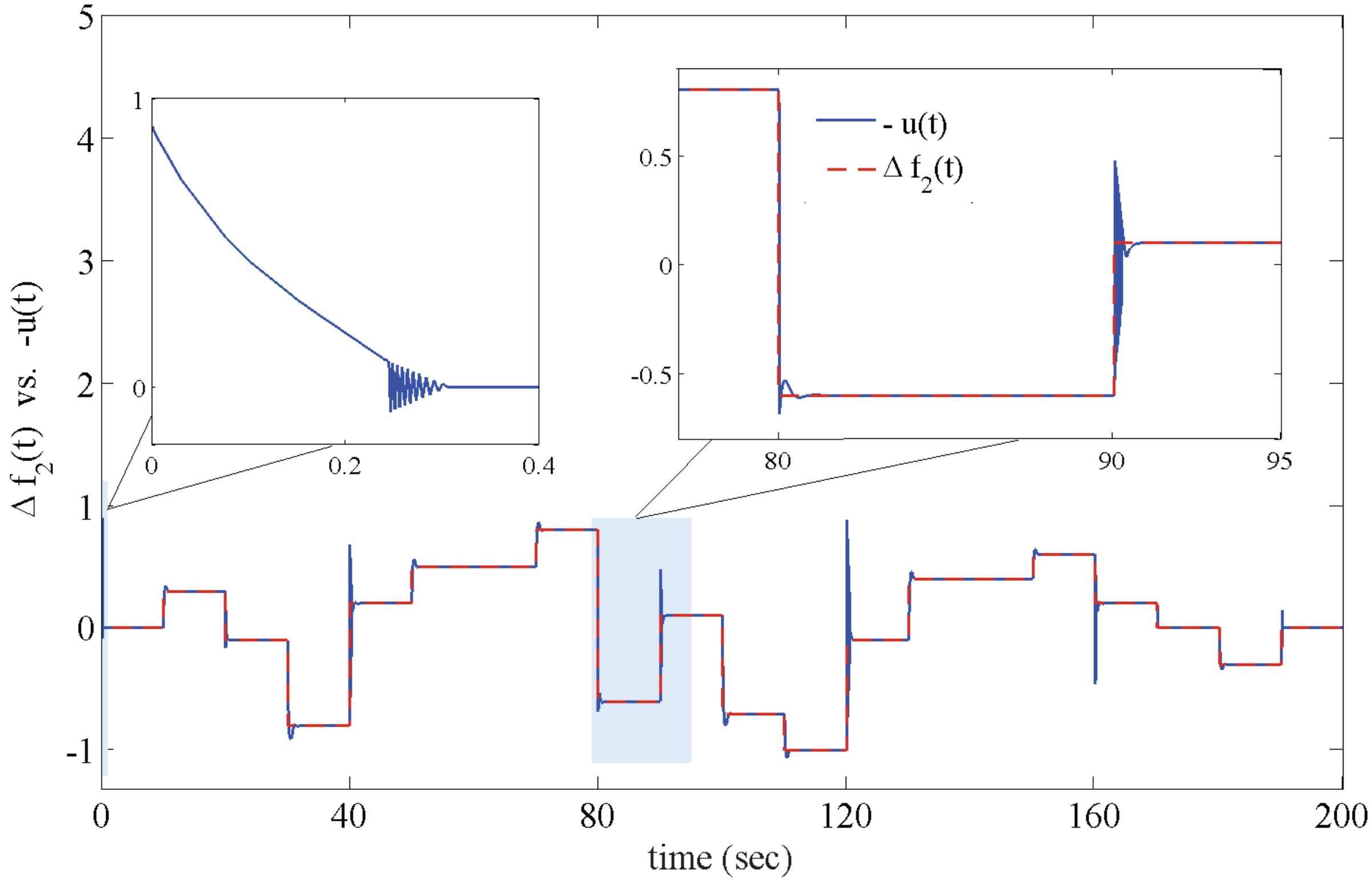}
		\caption{The adaptation performance for the square uncertainties $\Delta f_2(t)$}
		\label{fig:f1s}
	\end{figure}
	\begin{figure}[t]
		\centering
		\includegraphics[width = 2.4in]{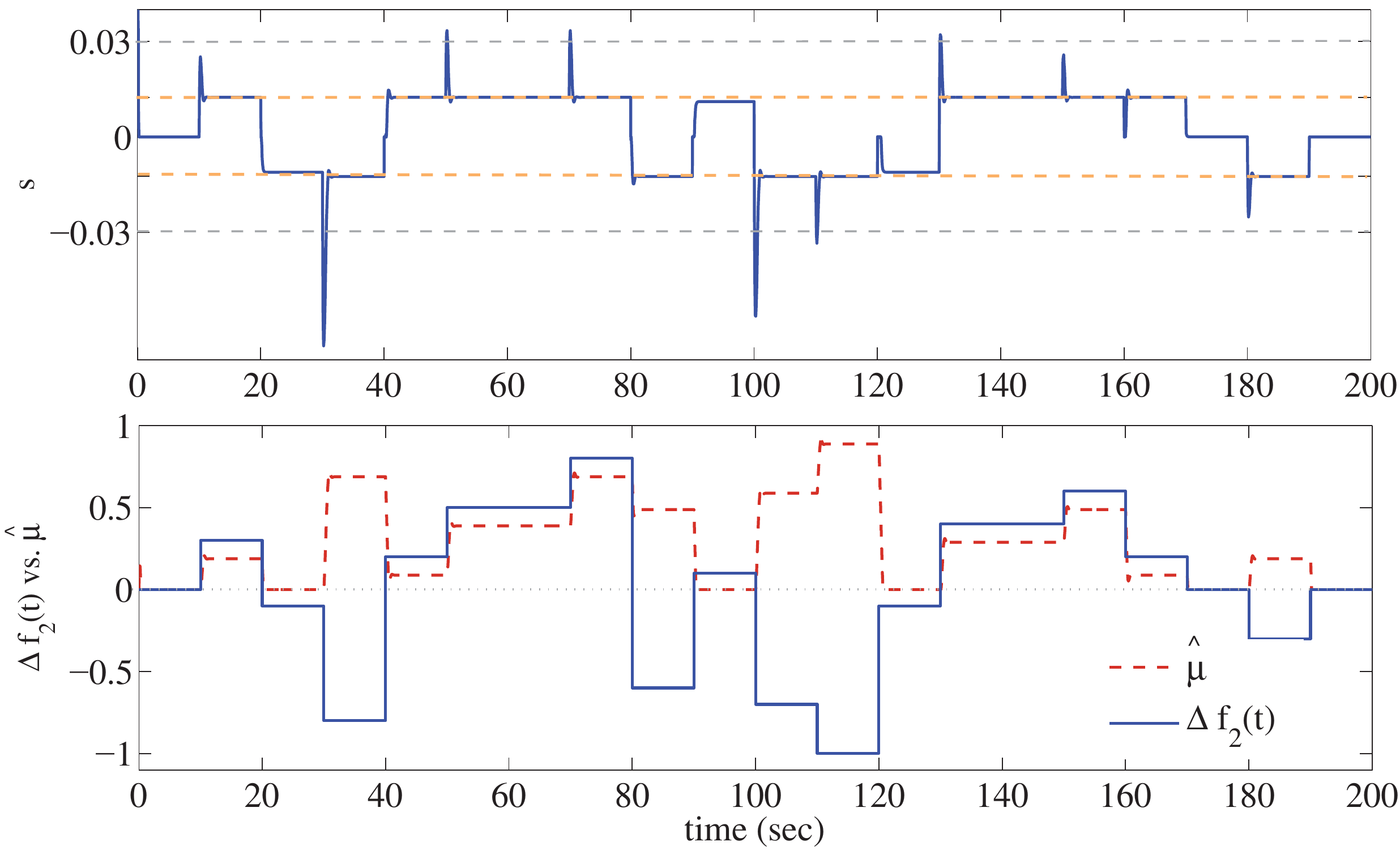}
		\caption{The convergence value of the sliding variables $s(t)$ vs. the adaptation gain $\hat{\mu}$ for the square uncertainties $\Delta f_2(t)$.}
		\label{fig:f2s}
	\end{figure}
	\begin{figure}[t]
		\centering
		\includegraphics[width = 2.4in]{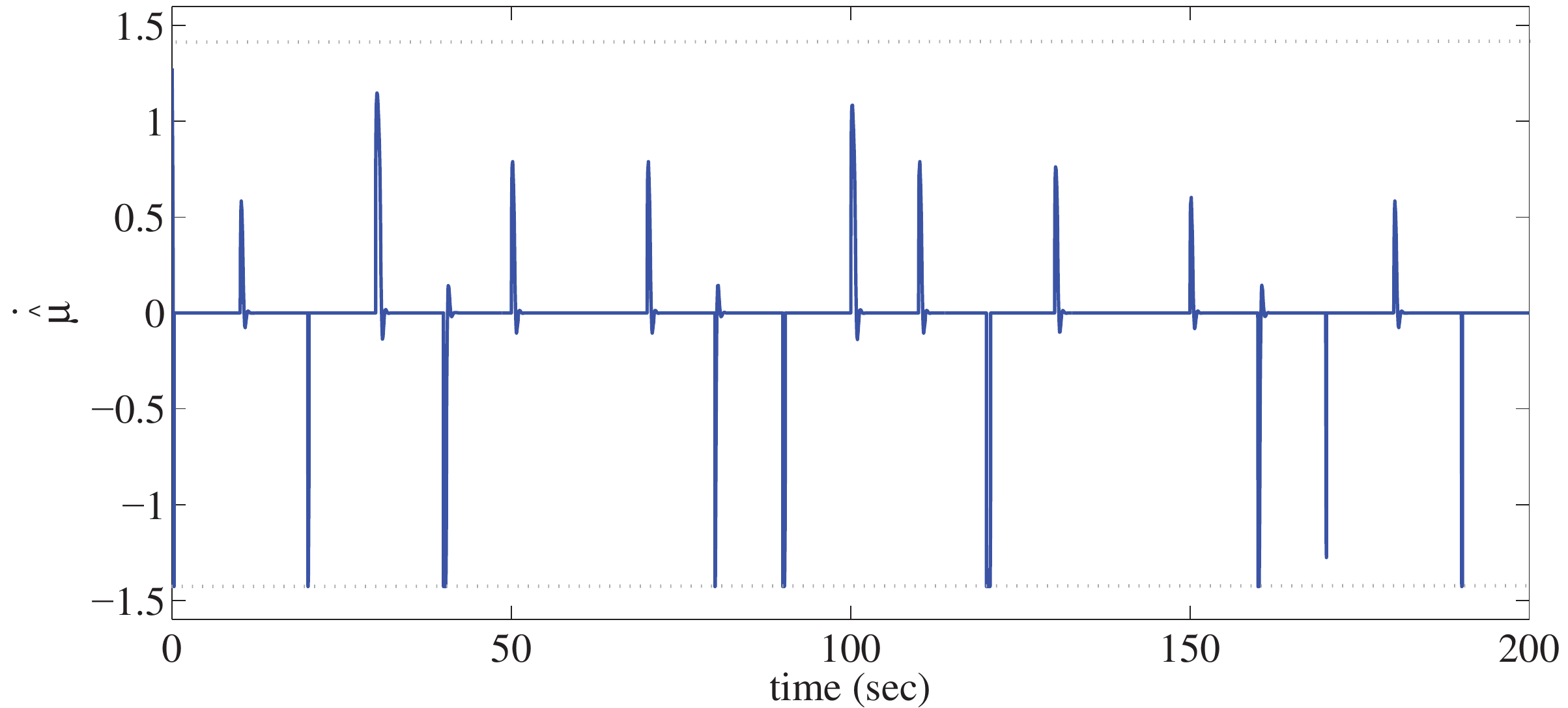}
		\caption{The changing rate of the adaptation gain, $\dot{\hat{\mu}}$ for the square uncertainties $\Delta f_2(t)$.}
		\label{fig:f3s}
	\end{figure}
	We obtain the time derivative of $V$ along the closed-loop system trajectories except $s = 0$ as 
	\begin{align}
	\dot{V}(s,\hat{\mu},t&)= \text{sgn}(s)\dot{s}_{\Delta}-\rho\dot{\hat{\mu}}(\mu-\hat{\mu}) \nonumber\\
	 &= (\Delta f\text{sgn}(s)-\hat{\mu}-k|s|)\Psi(s)-\rho \dot{\hat{\mu}}(\mu-\hat{\mu}) \label{eq:p}
	\end{align}
	where $\Psi(s)$ is defined by
	\begin{align}
	\Psi(s) = 1-\frac{2\phi^2}{(|s|+\phi)^2} \nonumber
	\end{align}
    for the sake of simplicity in later expressions.
    Fig. \ref{fig:delta_s} shows a plot of $\Psi(s)$ with varying thicknesses of the boundary layer. As we can see, the function $\Psi(s)$ intersects zero at the points $s = \pm \eta$. Moreover, it is positive when $s \in \mathcal{S}$ and negative outside. First, we consider the case of $\Psi(s) \ge 0$ which is $s \in \mathcal{S}$. $\dot{V}(s,\hat{\mu},t)$ becomes 
	\begin{align}
	\dot{V}(s,\hat{\mu},t) & \le (\mu-\hat{\mu}-k|s|)\Psi(s)-\rho\dot{\hat{\mu}}(\mu-\hat{\mu}) \nonumber\\
	& = (\mu-\hat{\mu})\left[\Psi(s)-\rho\dot{\hat{\mu}}\right]-k|s|\Psi(s). \label{eq:p1}
	\end{align}
	By setting $\dot{\hat{\mu}} = \Psi(s)/\rho$, we can eliminate the first term on the right hand side of equation (\ref{eq:p1}) and get the result of
	\begin{align}
		\dot{V}(s,\hat{\mu},t) \leq -k|s|\Psi(s) \leq 0 ~~ s \in \mathcal{S}. \label{eq:negative semi}
	\end{align}
	Notice that the result only gives us the update law of $\dot{\hat{\mu}} = \Psi(s)/\rho$ without the condition of $\hat{\mu}$ being non negative. If we substitute $\dot{\hat{\mu}}=0$ into (\ref{eq:p}), $\dot{V}$ will be indefinite. However, since $\dot{\hat{\mu}} \ge 0$ with $\hat{\mu}_0 >0$ for $s \in \mathcal{S}$, $\hat{\mu}$ can never be less than zero. In summary, we now have
	\begin{itemize}
		\item  $V(s,\hat{\mu})$ is monotonically increasing and bounded from below.
		\item  $\dot{V}(s,\hat{\mu},t) \le W(s,\hat{\mu}) = -k|s|\Psi(s) \le 0$ is negative semidefinite.
	\end{itemize} 
	for $(s,\hat{\mu}) \in \mathcal{S} \times \mathbb{R}^+$. 
	Note that Barbalat's Lemma is not applicable since we did not make any assumption on the uniformly continuity of the uncertainty.
	Alternatively, we can apply LaSalle's invariance principle (Theorem 2.2) from Barkana \cite{barkana2014defending} for the nonautonomous system. Based on satisfaction of assumption $1$ in Section 2.3 in \cite{barkana2014defending} for the boundedness of uncertainty (i.e. $|\Delta f| \le \mu$), we can conclude that all system trajectories are bounded and contained within the domain $\Omega = \{(s,\hat{\mu})\in \mathcal{S}| ~k|s|\Psi(s) = 0\}$ which implies
	\begin{align}
    (s,\hat{\mu}) \to (\pm\eta,~\mathbb{R}^+) \label{eq:theorem1result}
	\end{align}
	Now switch to the case of $\Psi(s) < 0$ for domain $\mathcal{S}'=\{s\in \mathbb{R},~|s|<\eta\}$. Substituting $\dot{\hat{\mu}} = \Psi(s)/\rho$ into (\ref{eq:p}), we have:
	\begin{align}
	\dot{V}(s,\hat{\mu},t) &= (\Delta f\text{sgn}(s)-k|s|-\mu)\Psi(s)\nonumber\\
	& \leq -(2\mu+k|s|)\Psi(s) \leq 2\mu+k|s| > 0\nonumber
	\end{align}
	where $\dot{V}$ is indefinite in the domain  $\mathcal{S}' \backslash \{0\}$ and undefined at $s=0$. Here, we cannot make any statement when $s \in \mathcal{S}'$. However, based on the result in (\ref{eq:theorem1result}), we can know that $s$ will approach $|s| = \eta$ when it is in $\mathcal{S}$.~~~~~~~~~~~~~~~~~~~~~~~~~~~~~\QED      
	\subsection{Convergence and Stability Analysis}
	In Section III.B, Theorem III.1, we only prove the convergence of $s$ to the boundary of $\mathcal{S}$ whenever $s \in \mathcal{S}$. However, there is no clear stability conclusion that can be drawn with respect to $\hat{\mu}$. Since $\Psi(s)$ is always positive in $\mathcal{S}$, it is possible that $\hat{\mu} \to \infty$ if $|s|$ never reaches the boundary within finite time. In this section, we show that $s(t)$ will reach the boundary of $\mathcal{S}$ with a finite $\hat{\mu}$ in finite time. Moreover, we can guarantee the trajectories of $s$ and $\hat{\mu}$ are bounded in steady state.

	\textbf{Theorem III.2:} Given the system (\ref{Problem}) implemented with (\ref{eq:control}) with initial conditions $(s_0\neq 0, \hat{\mu}_0>0)$ satisfying:
	\begin{align}
	|s_0|+\frac{1}{k} ~\hat{\mu}_0 = V'_0 > \frac{\sigma}{k}\text{, }~~~~\sigma = \mu + \frac{1}{k\rho}\label{eq:initial conditions}	
	\end{align}
	there exists a finite time $T$ such that
	\begin{align}
	~~~~~~~|s(t)|+\frac{1}{k}~\hat{\mu}(t) \le b,~~~~\forall t\ge T =\frac{1}{k} \text{ln}\frac{V'_0-\sigma/k}{b-\sigma/k} \nonumber
	\end{align}
	where $b$ is any number such that $\sigma /k < b < V'_0$.
	\begin{figure}[t]
		\centering
		\includegraphics[width = 2.8in]{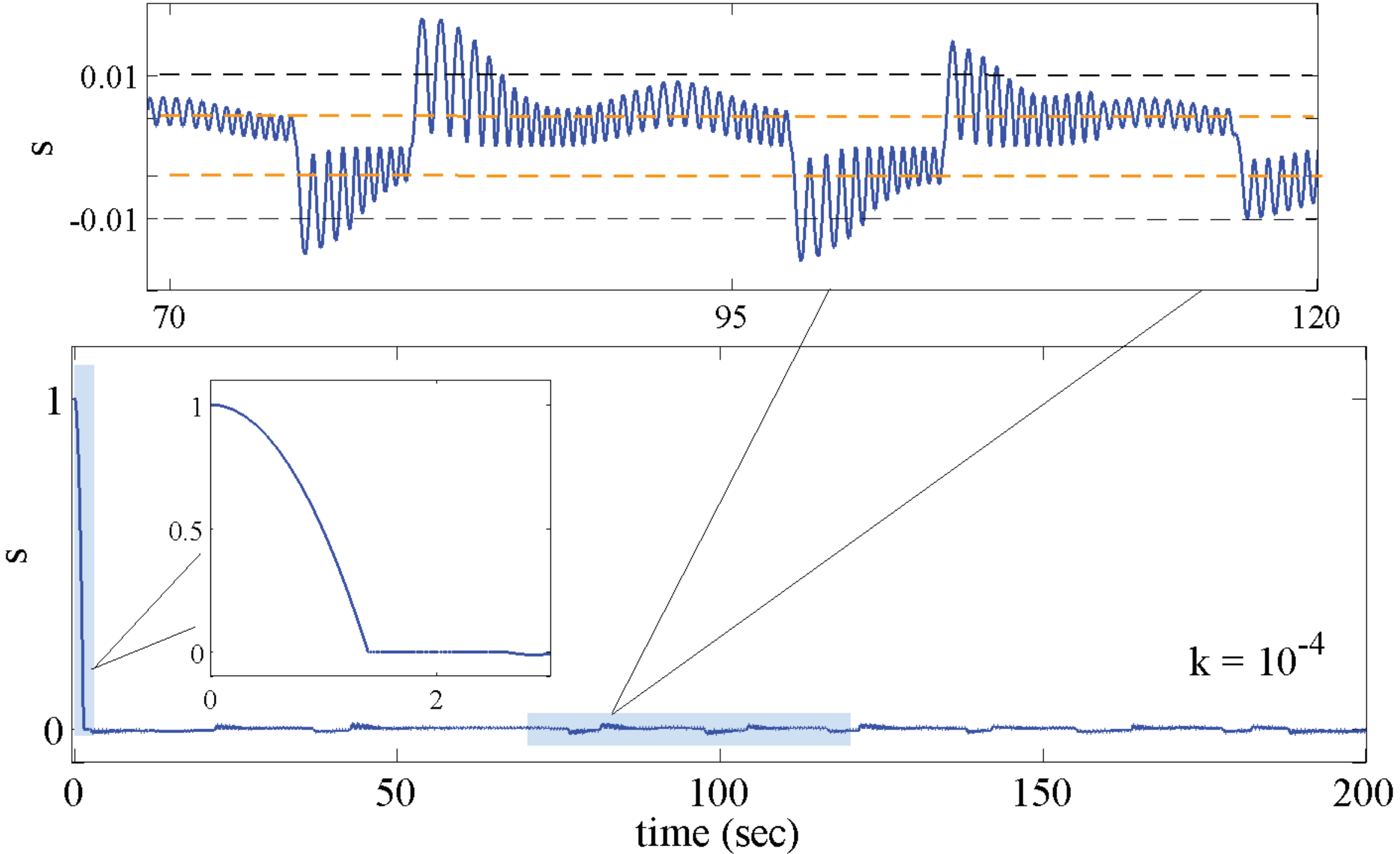}
		\caption{The sliding variable trajectories for the smooth uncertainties with a small feedback gain $k = 0.0001$}
		\label{fig:comf1}
	\end{figure}\\
	\textit{Proof:}
	Select another Lyapunov candidate:
	\begin{align}
	V'(t) = |s(t)|+\frac{1}{k}~\hat{\mu}(t) \nonumber
	\end{align}
	which is locally Lipschitz at $s=0$ and $\hat{\mu}=0$. Since $V'(t)$ is not differentiable everywhere, the upper right Dini derivative, $D^+{V}'(t)$, is introduced \cite{blanchini2009lyapunov}. By the assumptions of $\mathbf{l}(x,t)$ and $\mathbf{f}(x,t)$ in Section II.A, we know that the solution of equations (\ref{eq:sy}) and (\ref{eq:control}) exists and is absolutely continuous. Therefore, $D^+{V}'(t)$ is defined and the upper bound can be derived as:
	\begin{align}
		D^+{V}'(t) &= D^+\left[ |s(t)|+\frac{1}{k}\hat{\mu}(t)\right] \nonumber\\
		& \le |\Delta f| - \hat{\mu}-k|s| +\frac{1}{k\rho}\max\{\Psi(s),0\}\nonumber\\
		& \le \mu - \hat{\mu}-k|s| +\frac{1}{k\rho}\label{eq:kk}
	\end{align}
	according to the fact that $\Psi(s)$ is always bounded within the range of $[-1,1]$. Rewrite (\ref{eq:kk}) into $D^+V'(t) \le -kV'(t)+\sigma$. Then, the upper bound of the solution is given by
	\begin{align}
	& V'(t) \le e^{-kt}V'_0+\sigma \int\limits_0^te^{-k(t-\tau)}d\tau\nonumber\\
	\implies&|s(t)|+\frac{1}{k}~\hat{\mu}(t) \le b,~~~\forall t\ge T = \frac{1}{k}\text{ln}\frac{V'_0-\sigma/k}{b-\sigma/k}\label{eq:convergeProof}
	\end{align}
	where b is any number such that $\sigma /k < b < V'_0$.~~~~~~~~~~~~~~\QED
	
	\textit{\textbf{Remark:}} 
	As stated in Theorem III.1, $s$ is converging to the boundary $|s|=\eta$. From (\ref{eq:convergeProof}), we can know that the sum of $|s(t)|$ and $\hat{\mu}(t)$ is uniformly ultimately bounded \cite{khalil2002nonlinear} with ultimate bound $b$ after $T$. Therefore, we conclude that the $s(t)$ will reach the boundary of $\mathcal{S}$ with a finite $\hat{\mu}$ in finite time.
	Also, choosing the initial conditions satisfying (\ref{eq:initial conditions}) is not an issue in the implementation, since $\Psi(s)$ is always positive in $\mathcal{S}$ and the condition (\ref{eq:initial conditions}) will be met eventually for any initial setting of $s_0$ and $\hat{\mu}_0$. 
	
	\textbf{Theorem III.3:} Given the system (\ref{Problem}) implemented with the adaptation control law (\ref{eq:control}) the trajectories of $s$ are bounded within $|s(\mathbf{x},t)| < \delta$ after it first time reaches the domain $\mathcal{S}'$, where
	\begin{align}
	\delta = \sqrt{(2\eta)^2+\frac{\mu^2}{m}}-\eta \nonumber
	\end{align}
	and $m$ can be any value satisfying the following inequalities:
	\begin{align}
	m < \frac{\sqrt{2}}{\rho \phi}~~~~\text{and}~~~~\mu\sqrt{m}\leq\frac{1}{\rho}\Psi(\eta+\frac{\mu}{\sqrt{m}}). \label{eqn:m}
	\end{align}
		\textit{Proof:} According to the proof of Theorem III.1 and III.2, we get the result that $s$ will reach the boundary of $\mathcal{S}$ with a finite $\hat{\mu}$ in finite time.  
		To estimate the overshoot of $s$ after the first time it reaches the domain $\mathcal{S}'$, without loss of generality, consider the scenario when $s_0 = \eta^+$. Then, we choose an affine function to lower bound the original nonlinear adaptation law $\dot{\hat{\mu}} = \Psi(s)/\rho$ between the range $s = (\eta, \eta+\mu/\sqrt{m})$ and set $\hat{\mu}_0 = 0$, $k = 0$ in order to get the worst case response of $s$. 
		The system dynamics can be written as:
		\begin{align}
		\begin{cases}
		\dot{s}= -\hat{\mu}+\mu \\
		\dot{\hat{\mu}}= ms-m\eta
		\end{cases}. \label{eqn:upperbounddyn}
		\end{align}
		This yields
		\begin{align}
		s(t)&= (s_0+\eta)\cos(\sqrt{m}t)+\frac{\mu-\hat{\mu}_0}{\sqrt{m}}\sin(\sqrt{m}t)-\eta \nonumber\\
		&\leq \sqrt{(2\eta)^2+\frac{\mu^2}{m}}-\eta 
		\leq \eta+\frac{\mu}{\sqrt{m}}.
		\end{align}
		With this result, the requirements of $m$ in (\ref{eqn:m}) then are set. Because the slope of the adaptation law is equal to $\sqrt{2}/\rho\phi$ at $s = \eta$, the first requirement is set to allow the affine function to lower bound the nonlinear adaptation law between the range $s = (\eta, \eta+\mu/\sqrt{m})$. Then, the second requirement is set to ensure the validity of the dynamics ($\ref{eqn:upperbounddyn}$) within the range we claim.  ~~~~~~~~~~~~~~~~~~~~~~~~~~~~~~~~~~~~~~~~~~~~~~~~~\QED  
	\begin{figure}[t]
		\centering
		\includegraphics[width = 2.7in]{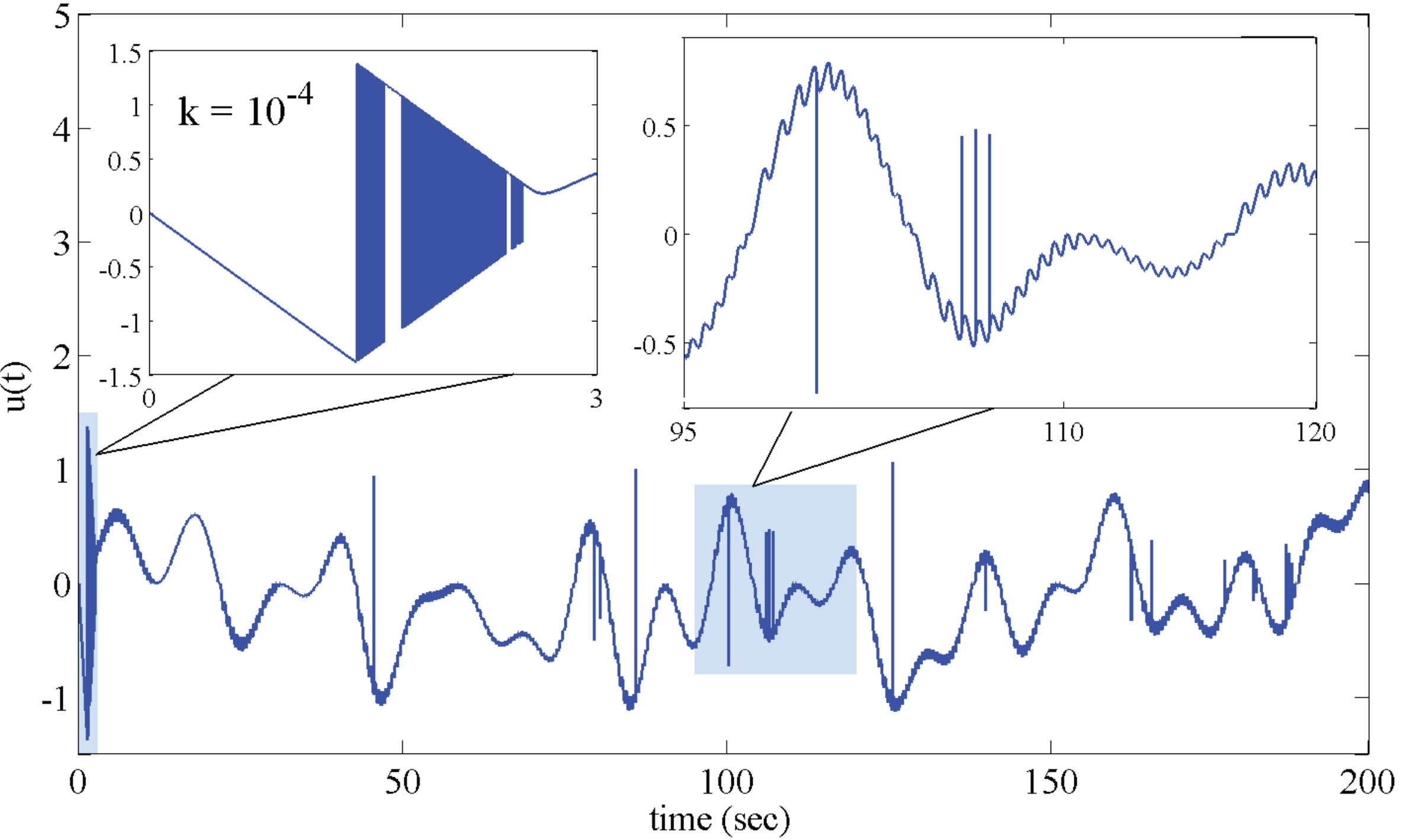}
		\caption{The control input for the smooth uncertainties with a small feedback gain $k = 0.0001$}
		\label{fig:comf2}
	\end{figure} 
	\section{Parameter Tuning and Implementation Issues}   
	In Section III, we have introduced three parameters $\phi$, $k$, $\rho$ for the new adaptation control strategy. Although the semi-global stability has been proven for all of them being positive in continuous-time, an adequate choice between each parameter is still needed for a good performance based on different scenarios. 
	\subsection{On the $\phi$-tuning}
    $\phi$ is a design parameter for the boundary of the domain $\mathcal{S}$ (i.e. $\pm\eta$) where the system trajectories will evolve during the steady state. Ideally, we would set $\phi$ as small as possible to have good tracking performance. However, setting $\phi$ too small will induce a large rate of change for the adaptation gain which may cause a high frequency chattering in both state responses and the control input. 
	We can observe that the slope of function $\Psi(s)$ around $\pm\eta$ becomes steeper as $\phi$ decreases in Fig. \ref{fig:delta_s}.  
	
	Another problem in the implementation would be the issue of discretization. Since nowadays many control algorithms are implemented using digital computers, an approximated discrete-time controller is commonly applied to the system. Setting $\phi$ too small may cause instability in the adaptation gain. To avoid this, we should follow the basic rule of thumb of allowing the system to sample roughly four times inside the domain $\mathcal{S}'$. As a result, the sampling rate will limit the parameter $\phi$ that we can choose. Therefore, we need to make a trade off between tracking performance and minimization of chattering through the choice of an adequate $\phi$.
	\begin{figure}[t]
		\centering
		\includegraphics[width = 2.4in]{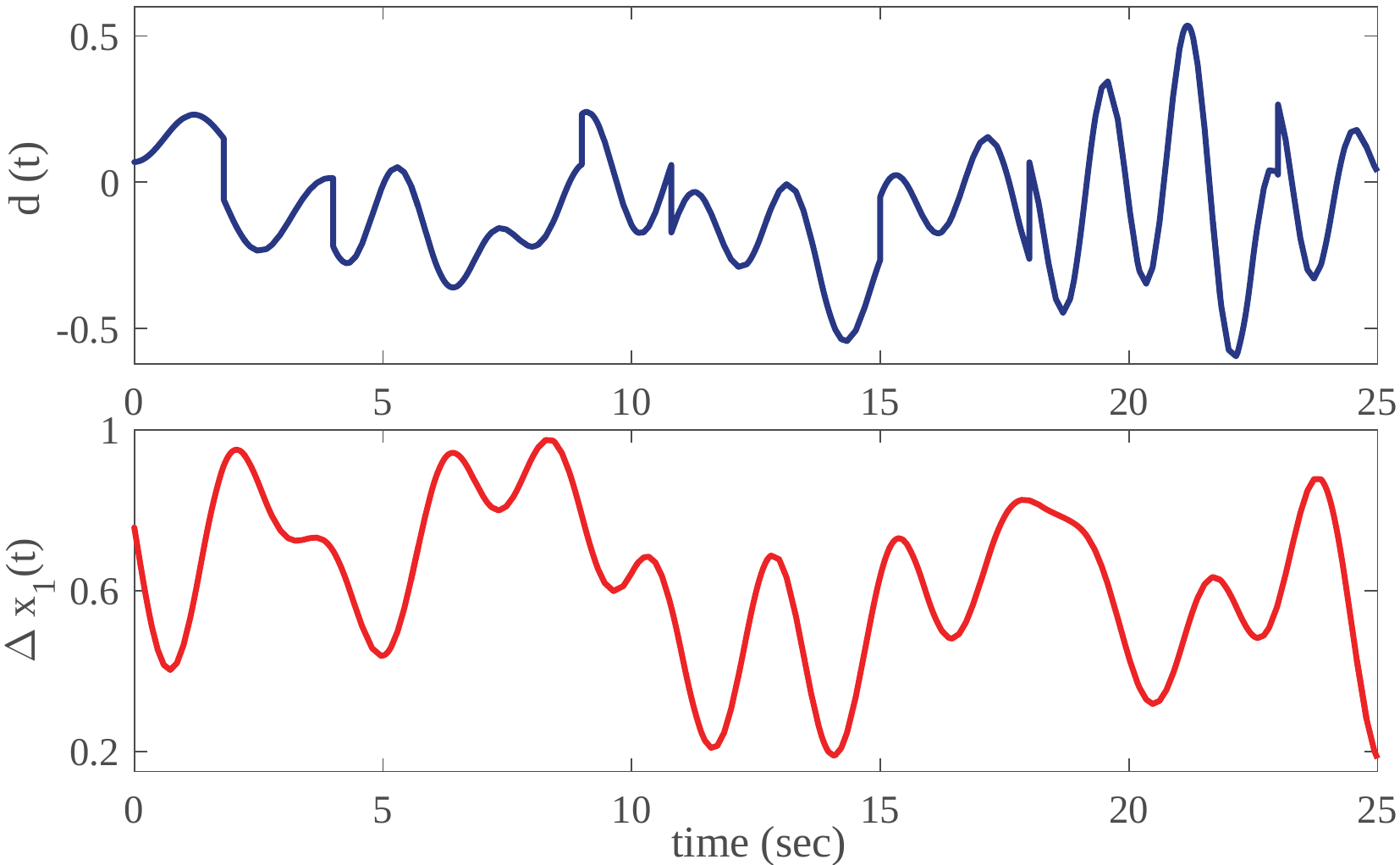}
		\caption{The multiplicative and additive  uncertainties, $\Delta x_1 (t)$ and $d(t)$.}
		\label{fig:disturbance}
	\end{figure}
	\begin{figure}[t]
		\centering
		\includegraphics[width = 2.5in]{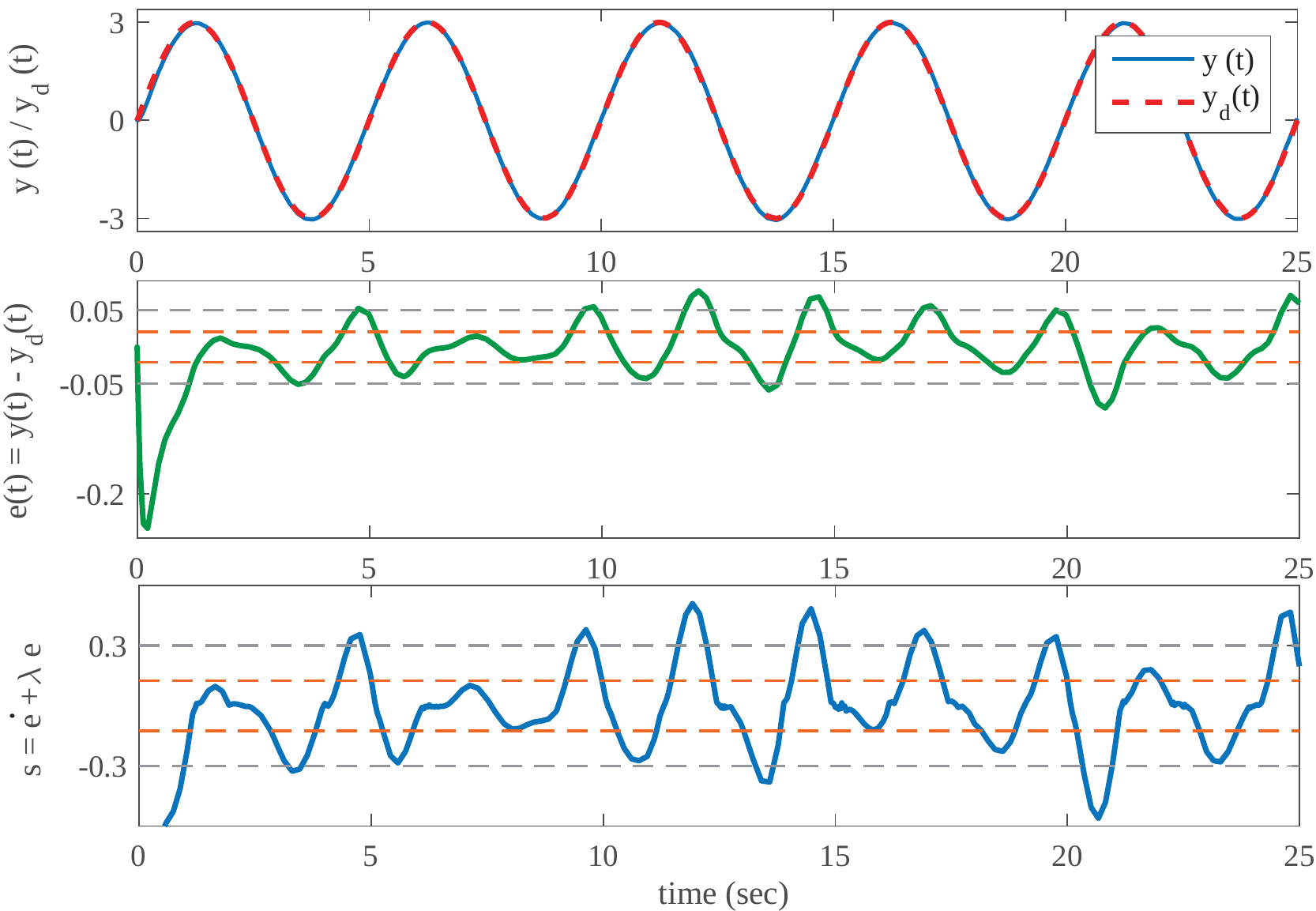}
		\caption{The tracking performance for the new proposed method.}
		\label{fig:tracking}
	\end{figure}
	\begin{figure}[t]
		\centering
		\includegraphics[width = 2.35in]{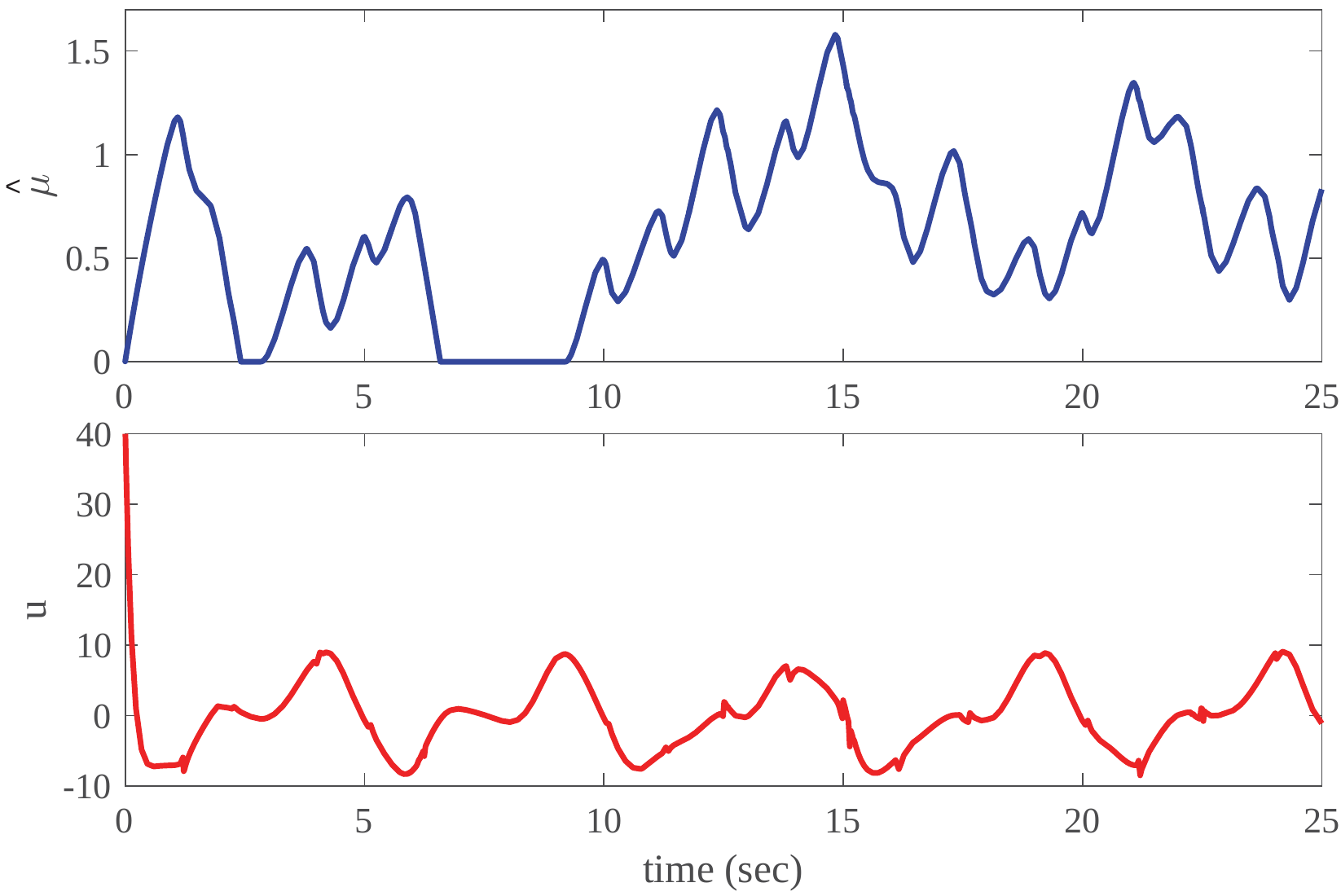}
		\caption{The adaptation gain $\hat{\mu}$ and the control input of the tracking problem for the new proposed method.}
		\label{fig:control}
	\end{figure}
	\begin{figure}[t]
		\centering
		\includegraphics[width = 2.4in]{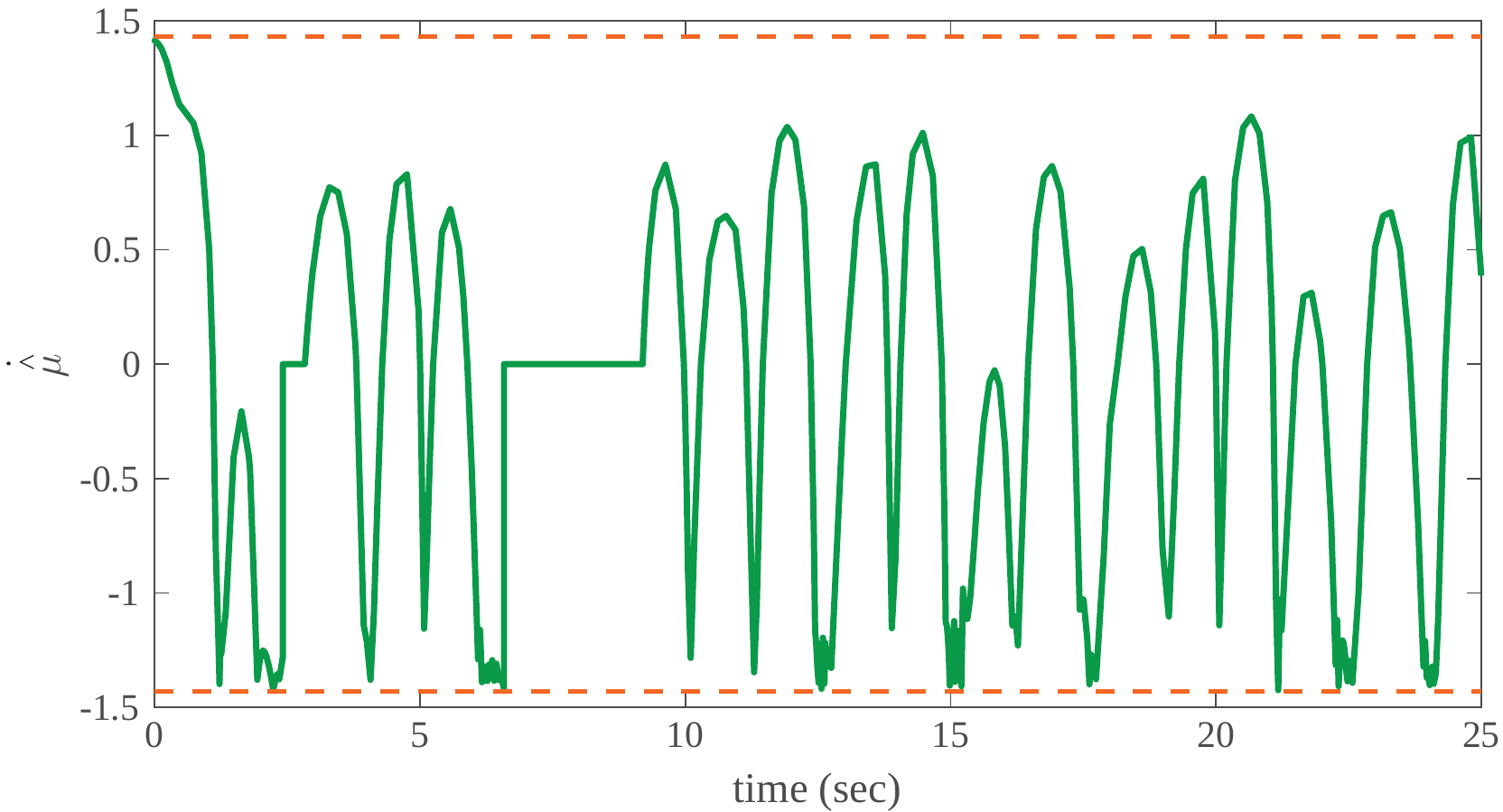}
		\caption{The changing rate of the adaptation gain, $\dot{\hat{\mu}}$, of the tracking problem for the new proposed method.}
		\label{fig:mudot}
	\end{figure} 
	\subsection{On the $k$-tuning}
	Compared with the standard sliding mode control law (\ref{eq:SMC control}), the new one described in (\ref{eq:control}) has an additional term $ks$, where $k$ is the design parameter for the feedback gain. 
	Having this additional control term benefits the overall performance since it will help speed up the convergence and smooth out the adaptation process.   
	Therefore, a higher value of $k$ ideally would be desired. However, in practice it should be limited by both actuator/unmodeled dynamics and the boundary thickness of $\mathcal{S}$.
	A high-gain control can easily excite unmodeled dynamics that could adversely affect the stability.
	For the second limitation of the boundary thickness, the condition of 
	\begin{align}
		0 \le k \le \frac{1}{\eta} \nonumber
	\end{align}  
	is required, since having too large of a feedback gain may lead to the system trajectories becoming confined inside an even smaller region of $\mathcal{S}$. Under this condition, the adaptation process will be terminated eventually as $\hat{\mu}$ goes to zero. 
	\subsection{On the $\rho$-tuning}	
	$\rho$ is the adaptation gain which is tuned based on the varying speed of the unknown uncertainties. Choosing a smaller $\rho$ allows a faster learning rate that can improve the adaptation process with high frequency uncertainties. However, we should notice that the smallest value of $\rho$ is limited by the actuation rate in application. 
	
	In conclusion, having smaller or larger values in both $\phi$ and $\rho$ or $k$ may be preferable, but all of them should be carefully chosen with suitable values to effectively avoid high control activity during the reaching phase and the adaptation process. 
	\section{SIMULATION}
	Two examples will be investigated in this section.
	First, we apply the adaptive control law given in (\ref{eq:control}) to a simple first-order system for a regulation problem in order to clearly demonstrate the properties of the new method.  
	Then, we again apply the control law to a higher order system with both multiplicative and additive uncertainties to a tracking problem. 
	\subsection{Regulation Problem}
	Consider the following system:
	\begin{align}
	\dot{x}=\Delta f(t) +u  \label{eq:ex1}
	\end{align}
    with $\Delta f(t)$ being bounded and unknown.
    Then, look for two different uncertainties (see in Fig. \ref{fig:df}) applied to this system: one is a smooth continuous function $\Delta f_1(t)$; the other is a sequence of square signals $\Delta f_2(t)$.  
    
    To implement the new adaptive control law, we first define the sliding surface as $s = x = 0$. Then, we choose the parameters $\phi = 0.01,~\rho = 1,~k=2,~\hat{\mu}_0=0.001,~x_0 = 1$ for the case with $\Delta f_1(t)$ uncertainty and $\phi = 0.03,~\rho = 0.7,~k=9,~\hat{\mu}_0=0.001,~x_0 = 0.1$ for the other case. Since the varying rates of the uncertainties are different, we choose a smaller $\rho$ in order to have a faster learning rate for the case of the square uncertainties. The effect is clearly seen. Fig. \ref{fig:f1} and Fig. \ref{fig:f1s} demonstrate the adaptation process works well under both low and high frequency uncertainties.
    The control input follows the external perturbation well particularly for slowly-varying uncertainties.
    \begin{figure}[t]
    	\centering
    	\includegraphics[width = 2.7in]{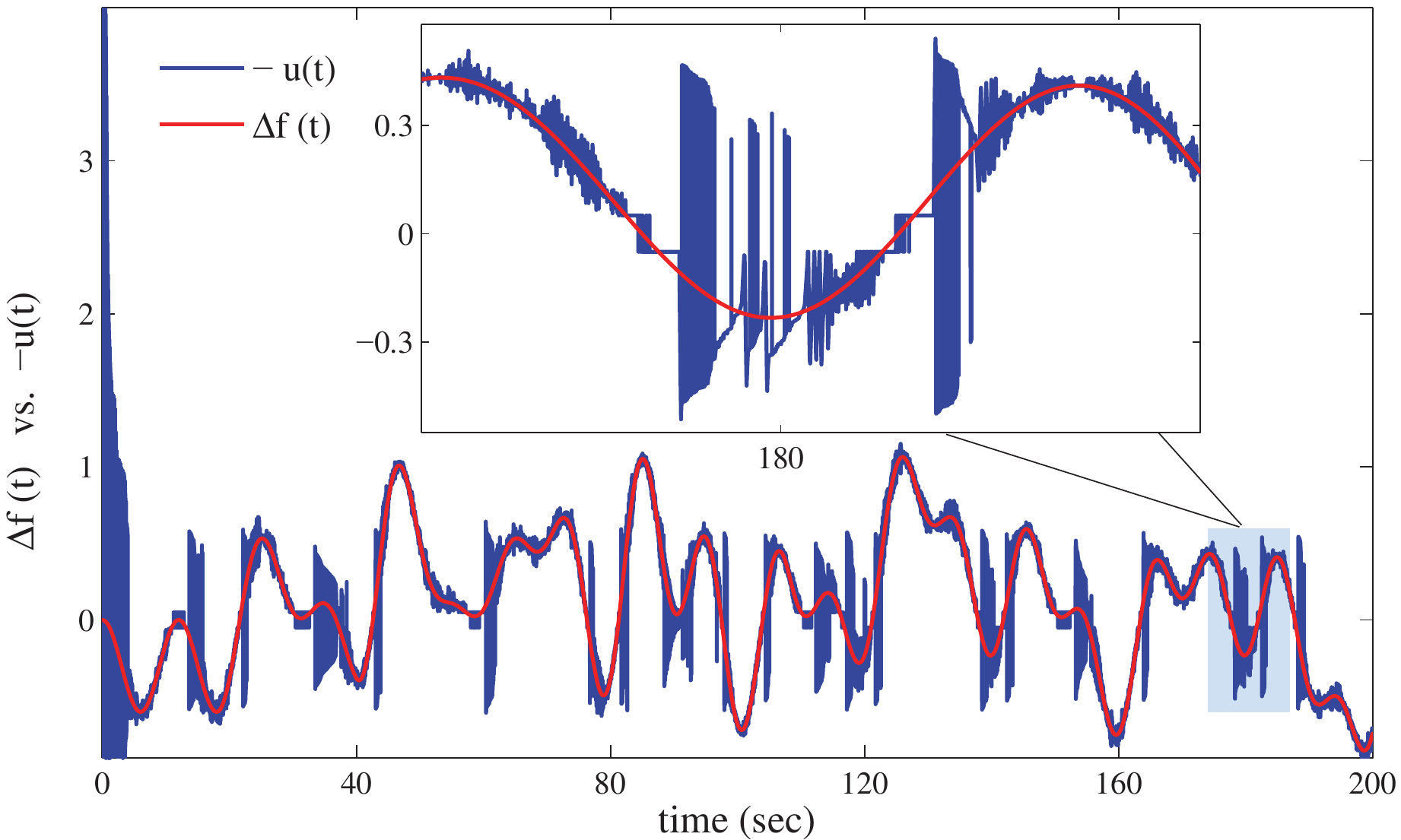}
    	\caption{The adaptation performance for the method proposed in \cite{plestan2010new}.}
    	\label{fig:old1}
    \end{figure}
    \begin{figure}[t]
    	\centering
    	\includegraphics[width = 2.7in]{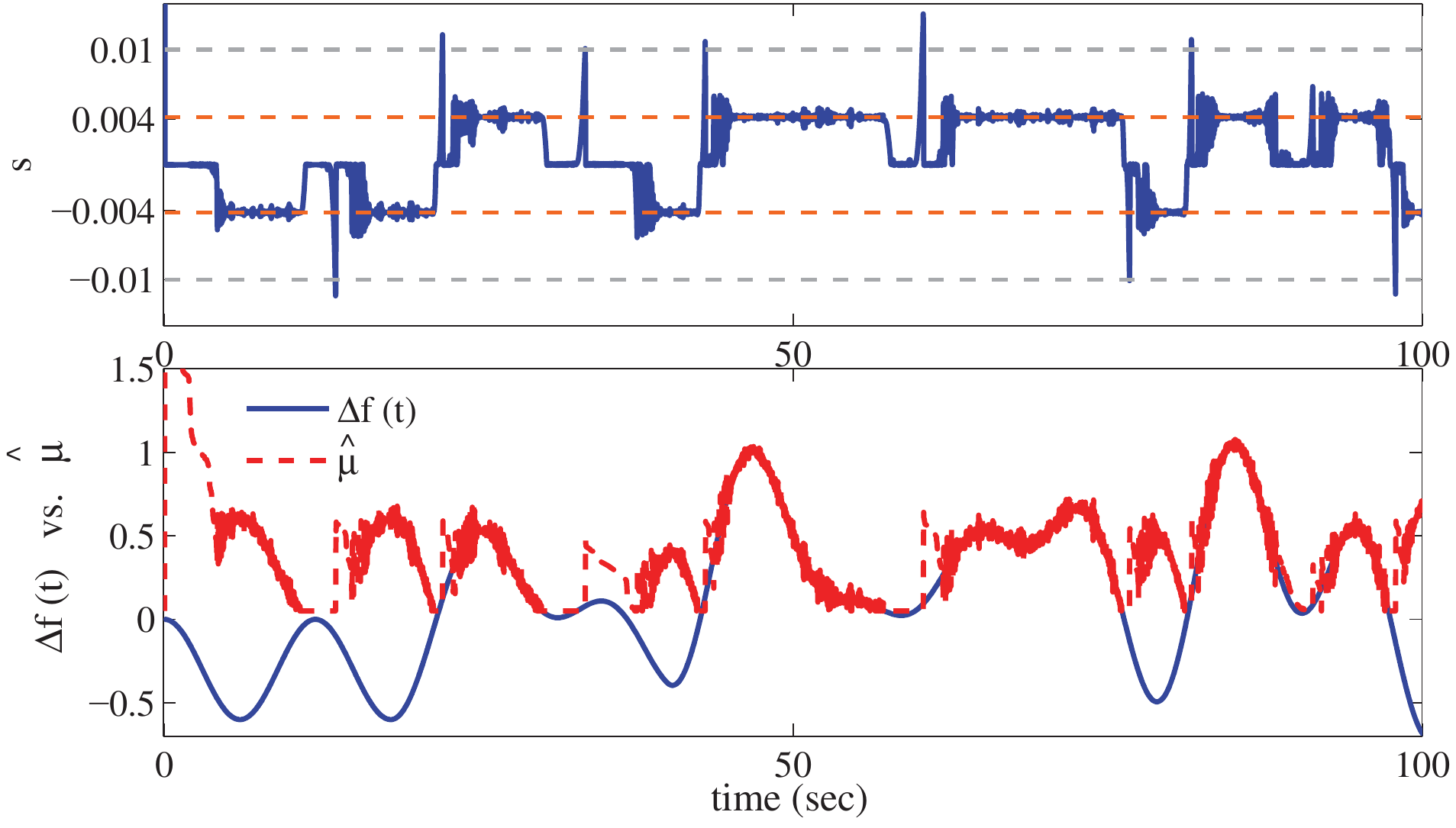}
    	\caption{The convergence value of the sliding variables $s(t)$ vs. the adaptation gain $\hat{\mu}$ for the method proposed in \cite{plestan2010new}.}
    	\label{fig:old2}
    \end{figure}
    \begin{figure}[t]
    	\centering
    	\includegraphics[width = 2.6in]{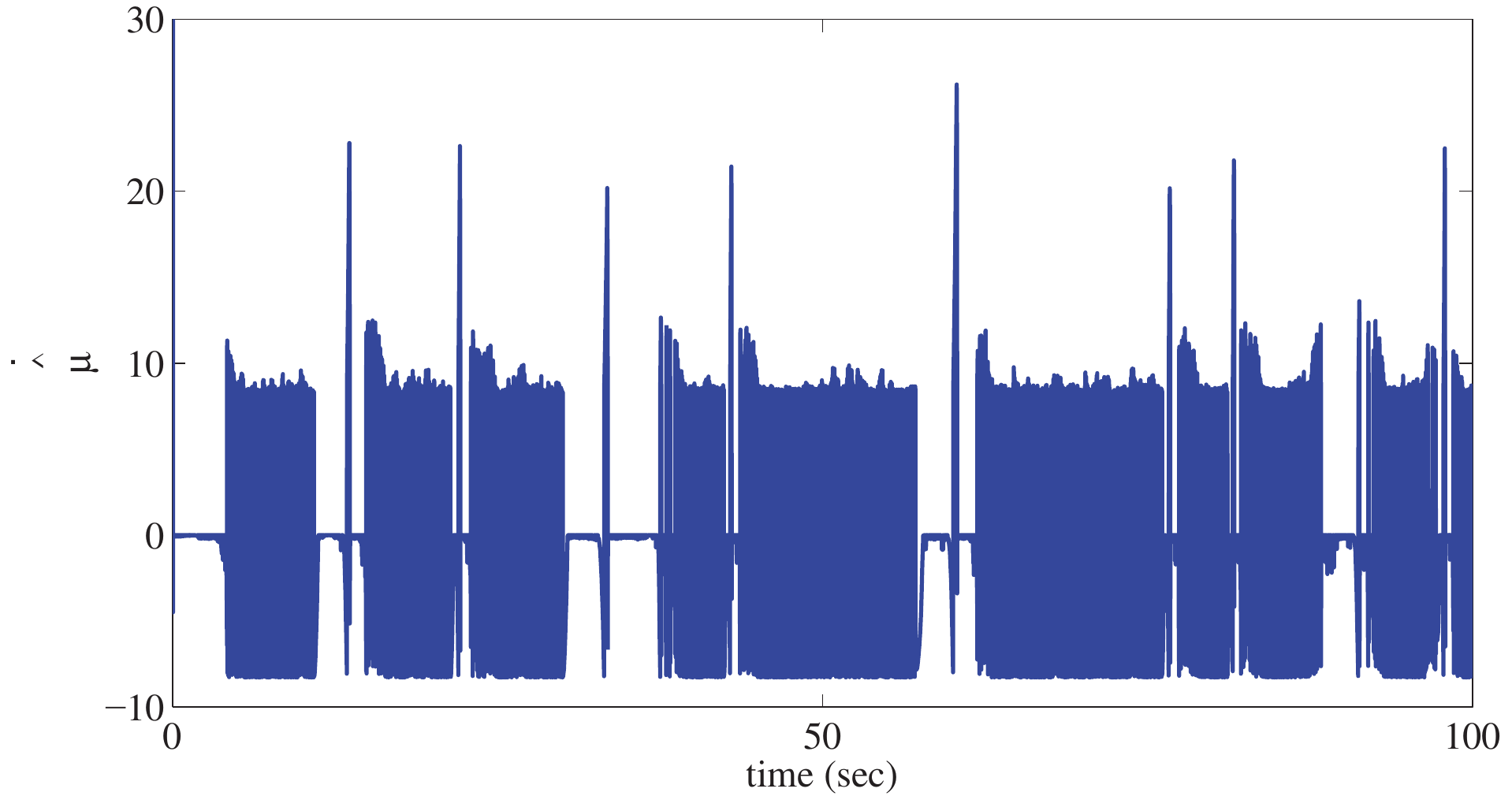}
    	\caption{The changing rate of the adaptation gain, $\dot{\hat{\mu}}$ for the method proposed in \cite{plestan2010new}.}
    	\label{fig:old3}
    \end{figure}
	Fig. \ref{fig:f2} and Fig. \ref{fig:f2s} show the response of state trajectories and the adaptation gain. We can see that in steady state, the sliding variable will evolve around the boundary of $\mathcal{S}$ (dotted orange line) without infinitely high frequency chattering. Since there are two convergence values $\pm (\sqrt{2}-1)\phi$ for the sliding variable, a connection between the convergence value of the sliding variable and the adaptation gain can be found in the figures. Because we restrict the adaptation gain to be always positive, the sliding variable will converge to the negative value of the boundary layer whenever the adaptation gain has the opposite sign of the current uncertainty.          
	
	The learning rates of the adaptation gain are shown in Fig. \ref{fig:f3} and Fig. \ref{fig:f3s}. We can notice that $\dot{\hat{\mu}}$ is always bounded within the region of $[-1/\rho, 1/\rho]$ since $|\Psi(s)|<1~\forall s\in \mathbb{R}$. Actually, it is one of the advantages of the proposed adaptation law compared with other methods. The main difference is that the update law is not a linear feedback law with respect to the sliding variable. The learning rate is limited and can be tuned by $\rho$. During the reaching phase, the state will converge with the auxiliary feedback term $ks$ instead of speeding up the increasing rate of the sliding gain.  
	This can help us smooth out the adaptation process and eliminate the oscillation behavior. We can verify the performance from Fig. \ref{fig:comf1} and Fig. \ref{fig:comf2} for the case of smooth uncertainties. The parameter setting is the same except $k = 0.0001$. The results show that both the state response and the control input signal perform worse when $k$ is small.  
	\subsection{Tracking Problem}
	Next, consider the following nonlinear system:
	\begin{align}
		\dot{x_1}&=x_2 \nonumber \\
		\dot{x_2} &=[x_1\Delta x_1(t)]x_2+\sin(x_1\Delta x_1(t))+d_1(t)+u\nonumber\\ 
		y& = x_1 \nonumber
	\end{align}
	with one multiplicative uncertainty, $\Delta x_1(t)$, and one additive uncertainty, $d(t)$, described in Fig. \ref{fig:disturbance}. The control objective is to apply the robust control law such that the output, $y$, tracks a reference signal, $y_{\text{d}} = 3\sin(0.4\pi\text{t})$. We first define $e = y-y_{\text{d}}$ and design a stable sliding surface as
	\begin{align}
	s = \dot{e}+\lambda e\text{, }~~~~~\lambda = 6 \nonumber
	\end{align}
    Then, apply the adaptive control law as
	\begin{eqnarray}
	\begin{aligned}
	u &= -x_1x_2-\sin x_1+\ddot{y}_{\text{d}}-\lambda( x_2-\dot{y}_{\text{d}})-ks-\hat{\mu}\text{sgn}(s) \\
	\dot{\hat{\mu}} &= \begin{cases}
	\frac{1}{\rho}\left[1-\frac{2 \phi^2}{(|s|+\phi)^2}\right]\text{~~if~~} \hat{\mu}\ge 0\\
	0\text{~~~~~~~~~~~~~~~~~~~~if~~} \hat{\mu} < 0 
	\end{cases} \hat{\mu}(0) = \hat{\mu}_0. \nonumber
	\end{aligned}
	\end{eqnarray}
    Assume the initial conditions of the states are all zero, $x_1=x_2=0$ and the parameters are set as $\rho = 0.7$, $k =5$, $\hat{\mu}_0 = 0.001$ and $\phi = 0.3$
	
    Fig. \ref{fig:tracking} demonstrates the tracking performance. As we can see in the second plot in Fig. \ref{fig:tracking}, the sliding variable, $s$, evolves around the boundary of $\mathcal{S}$ after it reaches the domain $\mathcal{S}'$. Additionally, since we have the result of
	\begin{align}
	s = 6 e \approx (\sqrt{2}-1)\phi \nonumber
	\end{align}
	in the steady state for $\dot{e}\approx 0$, we can know that the tracking error, $e$, will exhibit similar behavior as the sliding variable but with a scale of $1/6$. 
	Fig. \ref{fig:control} shows the simulation results of the adaptation gain and the control input. Although the overall uncertainty is unknown, we still can obtain a smooth adaptation process. Moreover, according to the analysis of the connection between the convergence value and the adaptation gain in the first example, we can even roughly reconstruct the overall uncertainty from the plots. The learning rate of the adaptation gain is shown in Fig. \ref{fig:mudot}. We can clearly see that $\dot{\hat{\mu}}$ is smooth  and always bounded within the range of $[-1/0.7,~1/0.7]$.        
	\section{Comparison}
	In this section, we compare the proposed adaptive sliding mode control with the one introduced in \cite{plestan2010new} by using the regulation example in Section V. Although the design concepts of these two methods are similar, the new one stands out for its smooth adaptation process without a high gain (i.e. $1/\rho$). We implement the adaptation law (\ref{eq:Adap}) in the case with the continuous uncertainty by setting the parameters $\bar{K} = 3000$ and $\kappa = 0.01$. Moreover, we choose $\epsilon = 0.01(\sqrt{2}-1)$ in order to have the same convergence standard for the comparison of these two methods. The simulation results in Fig.\ref{fig:old1}-\ref{fig:old3} display the closed-loop performance of the adaptive controller proposed in \cite{plestan2010new}. It appears that undesired chattering behaviors are introduced in both the control input and the sliding variable responses because of the large discontinuous switching rate in the adaptation gain (shown in Fig. \ref{fig:old3}). 
	
\begin{figure}[t]
		\centering
		\includegraphics[width = 2.8in]{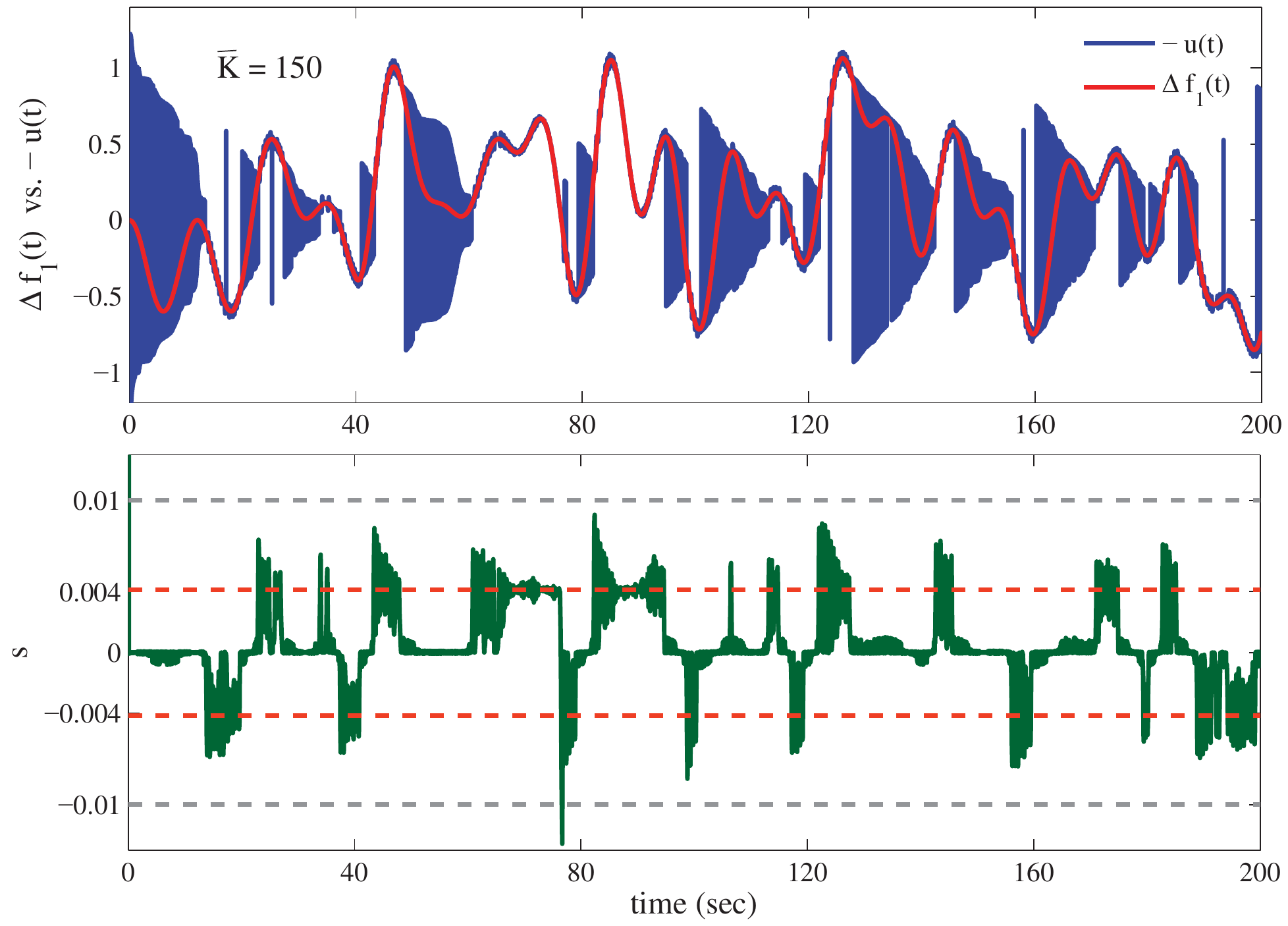}
		\caption{The performance of the control input canceling the uncertainty and the response of the $s(t)$ with a smaller $\bar{K}$.}
		\label{fig:old4}
	\end{figure}

\begin{figure}[t]
    	\centering
    	\includegraphics[width = 2.8in]{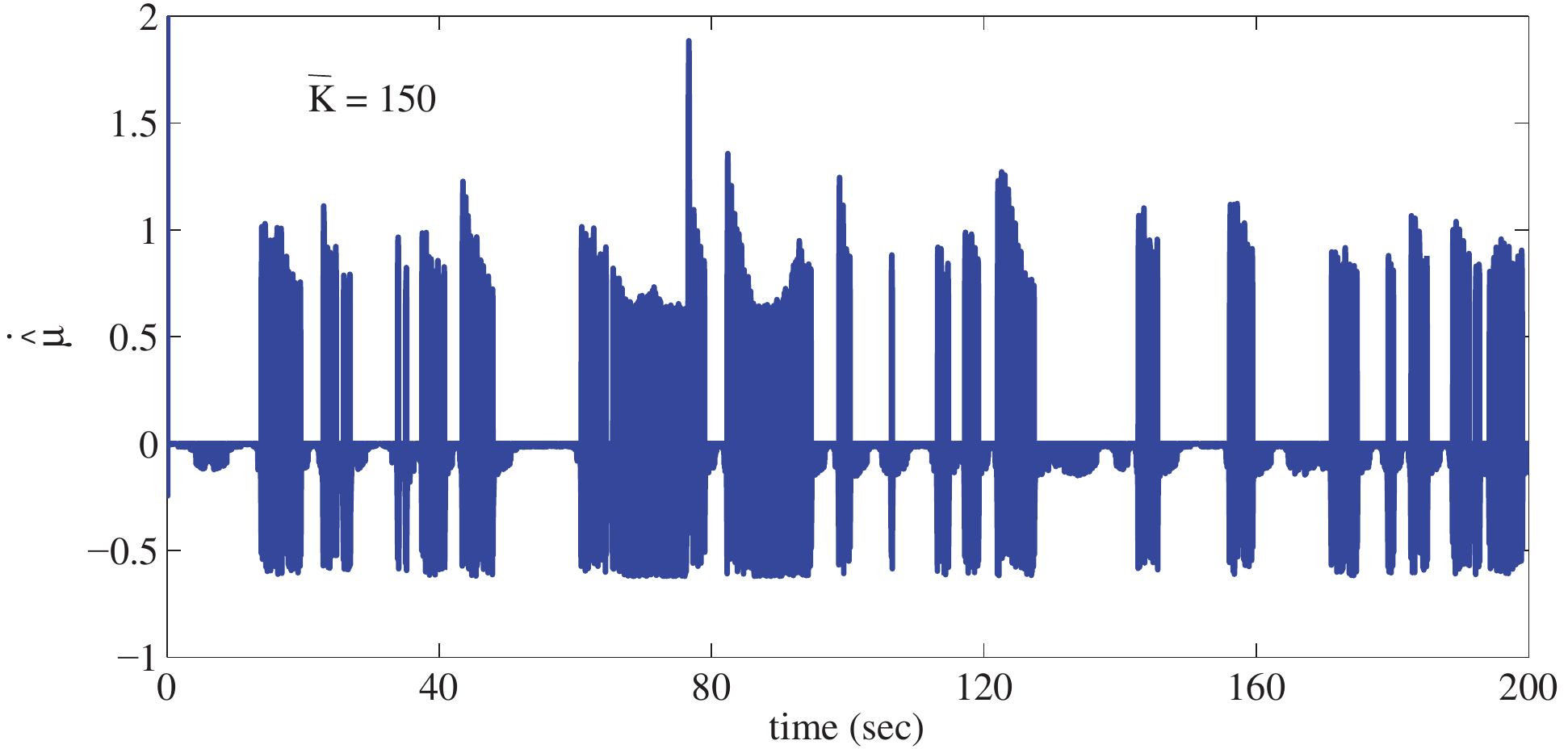}
    	\caption{The changing rate of the adaptation gain, $\dot{\hat{\mu}}$ for a smaller $\bar{K}$.}
    	\label{fig:old5}
    \end{figure}

	One suggestion for the alleviation of the chattering is to set the parameter $\bar{K}$ small, but as stated in Section III.A, the linear adaptation law with a smaller gain will make the decreasing rate even more insignificant inside the domain of $|s(\mathbf{x},t)|< \epsilon$, which fails to address the problem. Fig. \ref{fig:old4} and Fig. \ref{fig:old5} show the control effects with $\bar{K} = 150$. Although the changing rate of the adaptation gain becomes much smaller (Fig. \ref{fig:old5}), the chattering behaviors of the state and control input responses are not suppressed and even made worse with the small $\bar{K}$ (Fig. \ref{fig:old4}).                        
	
	\section{CONCLUSIONS}
	This paper proposed a new methodology of adaptive sliding mode control for a class of uncertain nonlinear systems. The algorithm utilizes the concept of the boundary layer. Based on the property that the system will hover inside and outside around the boundary region of $\mathcal{S}$, the adaptation law is designed such that the sliding gain will decrease and increase accordingly. Numerical examples illustrated the effect of the adaptation process. The process enables the determination of an adequate gain with respect to the current uncertainty. Semi-global stability of the closed-loop system with the adaptation gain is also guaranteed. Overall, this method achieves the minimum possible value of time-varying sliding mode control input and reduces the high-frequency chattering behavior without requiring knowing any knowledge of the uncertainties.
	
	%
	

\begin{thebibliography}{99}
	\bibitem{edwards1998sliding} C. Edwards and S. Spurgeon, Sliding mode control: theory and
	applications. Crc Press, 1998.
	\bibitem{utkin2009sliding} V. Utkin, J. Guldner, and J. Shi, Sliding mode control in electromechanical systems. CRC press, 2009, vol. 34.
	\bibitem{boiko2005analysis} I. Boiko, L. Fridman, and R. Iriarte, ``Analysis of chattering in continuous sliding mode control," in Proceedings of the 2005, American Control Conference, 2005. IEEE, 2005, pp. 2439-2444.
	\bibitem{yao1996smooth} B. Yao and M. Tomizuka, ``Smooth robust adaptive sliding mode control
	of manipulators with guaranteed transient performance," Journal of dynamic systems, measurement, and control, vol. 118, no. 4, pp. 764-775, 1996.
	\bibitem{utkin1992sliding} V. Utkin,``Sliding modes in control and optimization springer," New
	York, 1992.
	\bibitem{tseng2010chattering} M.-L. Tseng and M.-S. Chen,``Chattering reduction of sliding mode
	control by low-pass filtering the control signal," Asian Journal of control, vol. 12, no. 3, pp. 392-398, 2010.
	\bibitem{slotine1991applied} J.-J. E. Slotine, W. Li, et al., Applied nonlinear control. prentice-Hall Englewood Cliffs, NJ, 1991, vol. 199, no. 1.
	\bibitem{sastry2011adaptive} S. Sastry and M. Bodson, Adaptive control: stability, convergence and
	robustness. Courier Corporation, 2011.
	\bibitem{krstic1995nonlinear} M. Krstic, I. Kanellakopoulos, and P. V. Kokotovic, Nonlinear and
	adaptive control design. Wiley, 1995. 
	\bibitem{huang2008adaptive} Y.-J. Huang, T.-C. Kuo, and S.-H. Chang, ``Adaptive sliding-mode
	control for nonlinear systems with uncertain parameters," IEEE Transactions on Systems, Man, and Cybernetics, Part B (Cybernetics), vol. 38, no. 2, pp. 534-539, 2008.
	\bibitem{hussain2004adaptive} M. A. Hussain and P. Y. Ho, ``Adaptive sliding mode control with
	neural network based hybrid models," Journal of Process Control, vol. 14, no. 2, pp. 157-176, 2004. 
	\bibitem{hall2006sliding} C. E. Hall and Y. B. Shtessel, ``Sliding mode disturbance observer based
	control for a reusable launch vehicle," Journal of guidance, control, and dynamics, vol. 29, no. 6, pp. 1315-1328, 2006.
	\bibitem{utkin2013adaptive} V. I. Utkin and A. S. Poznyak, ``Adaptive sliding mode control with application to super-twist algorithm: Equivalent control method," Automatica, vol. 49, no. 1, pp. 39-47, 2013.
	\bibitem{dravzenovic1969invariance} B. Dra?zenovi´c, ``The invariance conditions in variable structure systems," Automatica, vol. 5, no. 3, pp. 287-295, 1969.
	\bibitem{khalil2002nonlinear} H. K. Khalil, ``Nonlinear systems," 2002.
	\bibitem{plestan2010new} F. Plestan, Y. Shtessel, V. Bregeault, and A. Poznyak, ``New methodologies
	for adaptive sliding mode control," International journal of control, vol. 83, no. 9, pp. 1907-1919, 2010.
	\bibitem{barkana2014defending} I. Barkana, ``Defending the beauty of the invariance principle," International Journal of Control, vol. 87, no. 1, pp. 186-206, 2014.
	\bibitem{blanchini2009lyapunov} F. Blanchini, ``Lyapunov methods in robustness. an introduction,"
	Lecture notes in Automatic Control, Bertinoro (Italy), 2009.	
	\end{thebibliography}

\end{document}